\documentclass{amsart}
\usepackage{amssymb, amsmath, amsfonts, amscd}

\input xypic

\theoremstyle{plain}
\newtheorem{theorem}{Theorem}[section]
\newtheorem{corollary}[theorem]{Corollary}
\newtheorem{lemma}[theorem]{Lemma}
\newtheorem{proposition}[theorem]{Proposition}

\newtheorem{example}[theorem]{Example}

\theoremstyle{definition}
\newtheorem{definition}[theorem]{Definition}

\theoremstyle{remark}

\numberwithin{equation}{theorem}

\newcommand{\Adef}{\mathcal{A}}
\newcommand{\K}{\operatorname{K}}

\renewcommand{\O}{\mathcal{O} }

\newcommand{\onab}{\overline{\nabla}}
\newcommand{\Iso}{\operatorname{Iso} }
\renewcommand{\mod}{\operatorname{Mod}}

\newcommand{\Hom}{\operatorname{Hom} }

\newcommand{\End}{\operatorname{End} }
\newcommand{\Ext}{\operatorname{Ext} }
\newcommand{\HH}{\operatorname{HH}}
\newcommand{\Tor}{\operatorname{Tor} }
\renewcommand{\H}{\operatorname{H} }
\newcommand{\Sym}{\operatorname{Sym} }

\newcommand{\Aut}{\operatorname{Aut}}

\newcommand{\nb}{\nabla}
\newcommand{\onb}{\overline{\nabla}}
\newcommand{\pb}{p_B}
\newcommand{\pl}{p_L}
\newcommand{\plf}{p_{L(f)}}

\newcommand{\Spec}{\operatorname{Spec} }

\newcommand{\C}{\operatorname{C} }
\newcommand{\Z}{\operatorname{Z}}
\newcommand{\B}{\operatorname{B}}
\newcommand{\D}{\operatorname{D}}
\newcommand{\Cx}{\mathbb{C} }
\newcommand{\Der}{\operatorname{Der} }

\renewcommand{\lg}{\mathfrak{g}}

\newcommand{\Q}{\mathbf{Q} }

\newcommand{\A}{R}
\renewcommand{\B}{A}
\newcommand{\U}{U^{ua}}
\newcommand{\conn}{\operatorname{Conn}}
\newcommand{\Bound}{\operatorname{B}}

\newcommand{\Mat}{\operatorname{Mat}}

\newcommand{\tL}{\tilde{L} }
\newcommand{\tp}{\tilde{p}}
\newcommand{\ts}{\tilde{s} }
\newcommand{\T}{\operatorname{T}}

\begin{document}

\title{Algebraic connections on projective modules with prescribed curvature}

\author{Helge Maakestad }
\email{$h\_maakestad@hotmail.com$}



\date{April 2014} 

\begin{abstract} In this paper we generalize some
results on universal enveloping algebras of Lie algebras to Lie-Rinehart algebras and twisted universal enveloping algebras
of Lie-Rinehart algebras. We construct for any Lie-Rinehart algebra $L$ and any 2-cocycle $f$ in $\Z^2(L,\B)$ 
the universal enveloping algebra $U(f)$ of type $f$. When $L$ is projective as left $\B$-module we prove a PBW-Theorem
for $U(f)$ generalizing classical PBW-Theorems. We then use this construction to give explicit constructions of a
class of finitely generated projective $\B$-modules with no flat algebraic connections.
One application of this is that for any Lie-Rinehart algebra $L$ which is projective
as left $\B$-module and any cohomology class $c$ in $\H^2(L,\B)$ there is a finite rank projective $\B$-module $E$ with
$c_1(E)=c$. 
Another application is to construct for any Lie-Rinehart algebra $L$ which is projective as left $\B$-module
a sub-ring $Char(L)$ of $\H^*(L,\B)$ - the characieristic ring of $L$. The ring $Char(L)$ ring is defined in terms of
the cohomology group $\H^2(L,\B)$ and has the property that it is a non-trivial subring of
the image of the Chern character $Ch_{\Q}:\K(L)_{\Q}\rightarrow \H^*(L,\B)$. 
We also give an explicit realization of the category of $L$-connections as a category of modules
on an associative algebra $\U(L)$.
\end{abstract}

\maketitle
\tableofcontents

\section{Introduction}

In the following paper we generalize classical notions on Lie algebras and universal enveloping algebras
of Lie algebras (see \cite{rinehart} and \cite{sridharan}) 
to Lie-Rinehart algebras and universal enveloping algebras of Lie-Rinehart algebras. We construct for any Lie-Rinehart 
algebra $L$
and any 2-cocycle $f$ in $\Z^2(L,\B)$, 
a generalized universal enveloping algebra $U(\B,L,f)$. The algebra $U(\B,L,f)$ is equipped with 
an increasing and decreasing filtration.
When $f=0$ we get Rinehart's universal enveloping algebra
$U(\B,L)$. We prove the following Poincare-Birkhoff-Witt Theorem for $U(B,L,f)$ (see Theorem \ref{pbw}): 

\begin{theorem} Let $L$ be projective as left $\B$-module.
There is a canonical isomorphism of graded $\B$-algebras
\[ \Sym^*_{\B}(L)\cong Gr(U(\B,L,f)) ,\]
where $Gr(U(\B,L,f))$ is the associated graded algebra of $U(\B,L,f)$ with respect to the ascending filtration.
\end{theorem}

As a consequence we get 
new examples of finitely generated projective modules with no flat algebraic connections. We also construct
families of (mutually non-isomorphic) finitely generated projective modules of arbitrary high rank using 
filtrations in the algebra $U(\B,L,f)$  (see Example \ref{family}).

Let $c\in \H^2(L,\B)$ be any cohomology class.
We use sub-quotients of the universal enveloping algebra  $U(\B,L,f)$ to  prove the following Theorem (see Theorem \ref{main}) 

\begin{theorem} Let $\A$ contain a field of characteristic zero and let $L$ be projective as left $\B$-module.
\label{main} There exist for every pair of integers $k,i\geq 1$ a finitely generated projective $\B$-module
$V(k,i)$ with the following property: 
\[ c_1(V(k,i))= c \in \H^2(L,\B).\]
Moreover $rk(V(k,i))=\binom{l+k+i-1}{l}-\binom{l+k-1}{l}$.
\end{theorem}

An application of this result is the following construction:
For any Lie-Rinehart algebra $L$ which is projective as left $\B$-module, 
there is a subring $Char(L)\subseteq \H^*(L,\B)$ which is defined in terms of the cohomology group $\H^2(L,\B)$.
The subring $Char(L)$ is a subring of the image $Im(Ch_{\Q})$ of the Chern character
\[ Ch_{\Q}:\K(L)_{\Q}\rightarrow \H^*(L,\B) .\]
The ring $Char(L)$ is non-trivial if and only if $\H^2(L,\B)\neq 0$. Hence 
we get a non-trivial subring of $Im(Ch_{\Q})$ whose definition does not involve
choosing generators of the Grothendieck group $\K(L)_{\Q}$. The problem of calculating generators of $\K(L)_{\Q}$ 
is an unsolved problem in general.

We also relate the cohomology group $\H^2(L,\B)$ where $(L,\alpha)$ is a Lie-Rinehart algebra which is projective as left
$\B$-module, to deformations of filtered associative algebras. Let $\Adef(\Sym_{\B}^*(L))$ be the deformation groupoid of the 
Lie-Rinehart algebra $(L,\alpha)$, parametrizing filtered associative algebras $(U,U_i)$ whose associated graded algebra
$Gr(U)$ is isomorphic to $\Sym_{\B}^*(L)$ as graded $\B$-algebra.
There is a one-to-one correspondence between $\H^2(L,\B)$ and the set
of isomorphism classes of objects in $\Adef(\Sym_{\B}^*(L))$ (see Theorem \ref{maindeform}).
As a Corollary it follows the category $\mod(U)$ of left $U$-modules is equivalent to the category of $L$-connections
of curvature type $f$, where $f$ is a 2-cocycle in $\Z^2(L,\B)$ (see Corollary \ref{modules}).
We also classify the morphisms in $\Adef(\Sym_{\B}^*(L))$ using the group $\Z^1(L,\B)$ 
(see Theorem \ref{morphisms}), hence the objects and morphisms of the 
deformation groupoid $\Adef(\Sym_{\B}^*(L))$ are determined by the cohomology group $\H^2(L,\B)$ and the group $\Z^1(L,\B)$.

Finally we realize in Theorem \ref{equiv} the category $\conn(L)$ of $L$-connections as a category of left 
modules on an associative ring
$\U(L)$. 

\begin{theorem} Let $L$ be a Lie-Rinehart algebra. There is an exact equivalence of categories
\[ \conn(L)\cong \mod(\U(L)) \]
preserving injective and projective objects.
\end{theorem}
The associative algebra $\U(L)$ is functorial in $L$. It is non-noetherian in general.

Using the algebra $\U(L)$ we give a definition of $\Ext$ and $\Tor$ modules
for arbitrary $L$-connections using the classical construction of $\Ext$ and $\Tor$ for modules on associative rings.
Cohomology and homology of connections was previously only defined for flat $L$-connections. The theory we 
define generalize Lie-Rinehart cohomology and homology for flat $L$-connections. 

We also classify the projective and injective objects in $\conn(L)$ using the algebra $\U(L)$. An $L$-connection 
$(W,\nabla)$ is a projective (resp. injective) object in $\conn(L)$ if and only if $W$ is a projective 
(resp. injective) $\U(L)$-module.
Hence projective objects in $\conn(L)$ are direct summands of free $\U(L)$-modules.

The algebra $\U(L)$ contains for every 2-cocycle $f$ for $L$ a 2-sided ideal 
$I_{L,f}$ with the property that there is an isomorphism
\[ \U(L)/I_{L,f}\cong U(\B,L,f).\]
Hence the algebra $\U(L)$ contains for every 2-cocycle $f$ an ascending and descending 
filtration with properties similar to the ascending and descending filtration of $U(\B,L,f)$. 
We moreover get a correspondence between 2-sided ideals in $\U(L)$ and curvature of $L$-connections (see Example 
\ref{twosided}). If $W$ is an $\U(L)$-module where $I_{L,f}$ annihilates $W$ it follows the corresponding connection
$\nabla$ has curvature type $f$ for some $f\in \Z^2(L,\B)$.
When $\B$ is noetherian and $L$ is finitely generated and projective, it follows $U(\B,L,f)$ is noetherian,
hence $\U(L)$ has many interesting non-trivial Noetherian quotients coming from geometry.

The main aim of the construction is to use it to construct moduli spaces of connections, cohomology of connections
and to apply this in the study of characteristic classes and the Chern-character.

\section{Lie-Rinehart cohomology and extensions}

In this section we extend well known results on Lie algebras, cohomology of Lie algebras and extensions to
cohomology of Lie-Rinehart algebras and extensions of Lie-Rinehart algebras. We give an interpretation of the 
cohomology groups $\H^i(L,W)$ for $i=1,2$ in terms of derivations of Lie-Rinehart algebras and equivalence
classes of extensions of Lie-Rinehart algebras.

Let in the following $h:\A \rightarrow \B$ be a map of commutative rings with unit. 

\begin{definition}
Let $L$ be a left $\B$-module and an $\A$-Lie algebra and let 
$\alpha:L\rightarrow \Der_{\A}(\B)$ be a map of left $\B$-modules and $\A$-Lie algebras.
The pair $(L, \alpha)$ is a \emph{Lie-Rinehart algebra} if the following equation holds
for all $x,y$ in $L$ and $a$ in  $\B$:
\[ [x,ay]=a[x,y]+\alpha(x)(a)y .\]
The map $\alpha$ is called the \emph{anchor map}.
Let $W$ be a left $\B$-module and let $\nabla: L\rightarrow \End_{\A}(W)$ be an $\B$-linear map. 
The map $\nabla$ is an \emph{L-connection} if the following equation holds
for all $x$ in $L$, $a$ in $\B$ and $w$ in $W$:
\[ \nabla(x)(aw)=a\nabla(x)(w)+\alpha(x)(a)w.\]
The category of $L$-connections is denoted $\conn(L)$.
\end{definition}

Let $(W,\nabla)$ be a connection.
Recall the definition
of the \emph{Lie-Rinehart complex} of the connection $\nabla$:
Let 
\[ \C^p(L, W):=\Hom_{\B}(\wedge^p_{\B} L, W) \]
with differentials
\[ d^p:\C^p(L,W) \rightarrow \C^{p+1}(L,W) \]
defined by
\[ d^p(\phi)(x_1\wedge \cdots \wedge X_p):=\sum_{k}(-1)^{k+1}\nabla(x_k)(\phi(x_1\wedge \cdots \wedge 
\overline{x_k}\wedge \cdots \wedge x_p))+\]
\[\sum_{i,j}(-1)^{i+j}\phi([x_i,x_j]\wedge x_1 \wedge \cdots \wedge \overline{x_i}\wedge \cdots \wedge 
\overline{x_j} \wedge \cdots \wedge x_p).\]

One checks the following:

\[ d^0(w)(x)=\nabla(x)(w),\]
\[ d^1(\phi)(x\wedge y)=\nabla(x)(\phi(y))-\nabla(y)(\phi(x))-\phi([x,y]),\]
and
\[ d^1(d^0(w))(x\wedge y)=R_\nabla(x\wedge y)(w), \]
where
\[ R_\nabla(x\wedge y)=[\nabla(x),\nabla(y)]-\nabla([x,y]).\]

\begin{definition} Let $R_{\nabla}$ be the \emph{curvature} of the connection $\nabla$.
We say $\nabla$ is \emph{flat} is $R_{\nabla}$ is zero.
\end{definition}

One checks that the sequence  of groups and maps given by 
$( \C^p(L,W), d^p)$ is a complex of $\A$-modules if and only if the curvature $R_\nabla$  is zero. 

\begin{definition} Let $(W,\nabla)$ be a flat $L$-connection. Let $\Z^i(L,W):=ker(d^i)$ and 
$\Bound^i(L,W):= im(d^{i-1})$.
 Let for all $i\geq 0$ $\H^i(L,W):=\Z^i(L,W)/\Bound^i(L,W)$ be the \emph{i'th Lie-Rinehart cohomology group}
of $L$ with values in $(W,\nabla)$
\end{definition}

It follows the abelian group $\H^i(L,W)$ is a left $\A$-module.

In his PhD-thesis \cite{rinehart} Rinehart introduced the universal enveloping algebra $U(\B,L)$ for a Lie-Rinehart
algebra $L$ and proved a PBW-Theorem for $U(\B,L)$ in the case when $L$ is a projective $\B$-module. He also proved various
general results on the cohomology groups $\H^i(L,W)$ using the algebra $U(\B,L)$. This was the first systematic
study of the algebra $U(\B,L)$ and the cohomology groups $\H^i(L,W)$, hence the name  \emph{Lie-Rinehart cohomology}.

The complex $\C^p(L,W)$ has many names: The Lie-Rinehart complex, the Chevalley-Hochschild
complex, the Lie-Cartan complex, the Chevalley-Eilenberg complex etc. It was known prior to Rineharts paper \cite{rinehart}
that for a real smooth manifold $M$ with $\O(M)$ the algebra of real valued smooth functions and 
$\mathcal{L}=\Der_{\mathbf{R}}(\O(M))$ the Lie-algebra of derivations of $\O(M)$, it follows there is an isomorphism
\[ \H^i(\mathcal{L},\O(M)) \cong \H^i_{sing}(M,\mathbf{R}) \]
for all $i\geq 0$, where $\H^i_{sing}(M,\mathbf{R})$ is singular cohomology of $M$ with real coefficients.
I have not found a reference to this result but it appears this was considered ``folklore'' at the time.
The result also follows from results in Rineharts paper \cite{rinehart} and the classical DeRham Theorem.

If $A$ is a regular algebra of finite type over the complex numbers and $X:=\Spec(A)$ the affine scheme associated 
to $A$ one may consider $X(\Cx)$ - the underlying complex algebraic manifold of $X$. It follows from \cite{grothendieck}
there is an isomorphism
\[ \psi_i:\H^i(\Der_{\Cx}(A),A)\cong \H^i_{sing}(X(\Cx),\Cx) \]
for all $i\geq 0$, 
where $\H^i_{sing}(X(\Cx),\Cx)$ is singular cohomology of $X(\Cx)$ with complex coefficients. The existence of the 
isomorphisms $\psi_i$ 
was proved for a general smooth affine algebraic variety by Grothendieck in \cite{grothendieck}. It was known
to exist for affine homogeneous spaces by Hochschild and Kostant. It is remarkable that the complex $\C^*(\Der_{\Cx}(A),A)$
which is a purely algebraic object, calculates singular cohomology of $X(\Cx)$. 

For a field $k$ of characteristic $p>0$ one may use the groups $\H^i$ to construct a $p$-adic cohomology
theory of varieties over $k$  with properties similar to crystalline cohomology 
(see \cite{berthelot} for an introduction to crystalline cohomology). Lie-Rinehart cohomology
is known to generalize several other cohomology theories: Algebraic De Rham cohomology, 
logarithmic De Rham cohomology, poisson cohomology and  Lie algebra cohomology.
It is well known Lie-Rinehart cohomology does not calculate singular cohomology of a 
manifold $M$ with integer or rational coefficients.

In this section we are interested in the group $\H^i(L,W)$ for $i=1,2$ where $(W,\nabla)$ is a flat connection. 

We get a map
\[ d^2:\C^2(L,W)\rightarrow \C^3(L,W)\]
where for any element
\[ f\in \C^2(L,A):=\Hom_{\B}(\wedge^2_{\B} L,W) \]
it follows 
\[d^2(f)(x_1 \wedge x_2 \wedge x_3)=\nabla(x_1)(f(x_2\wedge x_3))-\nabla(x_2)(f(x_1\wedge x_3)) 
+\nabla(x_3)(f(x_1 \wedge x_2)) \]
\[ -f([x_1,x_2]\wedge x_3)+f([x_1,x_3]\wedge x_2)-f([x_2,x_3]\wedge x_1). \]
It follows $\Z^2(L,W)$ is the set of $\B$-bilinear maps
\[ f:L\times L \rightarrow W \]
satisfying $f(x,x)=0$ for all $x\in L$ and such that $d^2(f)=0$. 

\begin{definition}
Let $\alpha:L\rightarrow \Der_{\A}(\B)$ and $\tilde{\alpha}:\tilde{L}\rightarrow \Der_{\A}(\B)$
be Lie-Rinehart algebras. 
Let
\[ p:\tilde{L}\rightarrow L \]
be a map of left $\B$-modules and $\A$-Lie algebras.
We say $p$ is a map of \emph{Lie-Rinehart algebras} if $\alpha \circ p =\tilde{\alpha}$.
\end{definition}

Let $p:\tilde{L} \rightarrow L$ be a surjective map of Lie-Rinehart algebras and let $W=ker(p)$. It follows $W$ is a
sub-$\B$-module and a sub-$\A$-Lie algebra of $\tilde{L}$. We may view $W$ as a \emph{trivial Lie-Rinehart algebra}
with trivial Lie-product and trivial anchor map. Such a Lie-Rinehart algebra is sometimes called a
\emph{totally intransitive Lie-Rinehart algebra}.
We get an \emph{exact sequence} of Lie-Rinehart algebras 
\[ W \xrightarrow{i} \tilde{L}\xrightarrow{p} L  .\]
This means the map $i$ is an injection, $p$ a surjection and $Im(i)=Ker(p)$.
Moreover $i$ and $p$ are maps of left $\B$-modules and $\A$-Lie algebras.
Define the following action:
\begin{align} 
&\label{flat} \tilde{\nabla}:\tilde{L}\rightarrow \End(W) 
\end{align}
by
\[ \tilde{\nabla}(z)(w):=[z,w]\]
where $[,]$ is the Lie-product on $\tilde{L}$ and $z\in \tilde{L}, w\in W$.
It follows  the map $\tilde{\nabla}$ is a flat $\tilde{L}$-connection on $W$.
Assume $W=ker(p) \subseteq \tilde{L}$ is an abelian sub-algebra of $\tilde{L}$. 
Assume $z\in \tilde{L}$ is an element with $p(z)=x\in L$. Let $w\in W$. Define the following map:

\begin{align}
&\label{flat2} \rho:L\rightarrow \End(W) 
\end{align}
by
\[ \rho(x)(w):=[z,w].\]
Assume $p(z')=x$. It follows $z'=z+v$ where $v\in W$. We get $[z+v,w]=[z,w]+[v,w]=[z,w]$. Hence the 
element $\rho(x)\in \End(W)$ does not depend on choice of the element $z$ mapping to $x$. It follows $\rho$ is a well
defined map. One checks that $\rho$ is a $\B$-linear map
\[ \rho:L\rightarrow \End_{\A}(W).\]
One checks the map $\rho$ is a flat $L$-connection $W$.

\begin{definition}
Fix a flat connection
\[ \nabla:L\rightarrow \End_{\A}(W) \]
on the Lie-Rinehart algebra $L$ and assume $p:\tilde{L}\rightarrow L$ is a surjective map of Lie-Rinehart algebras.
Assume $W=ker(p)$ is an abelian sub-algebra of $\tilde{L}$. Assume the induced connection
\[ \rho:L\rightarrow \End_{\A}(W) \]
 from (\ref{flat2}) equals $\nabla$.
The extension
\[  W \rightarrow \tilde{L}\rightarrow L \]
is an \emph{extension of $L$ by the flat connection $(W,\nabla)$}.
Two extensions $L_1,L_2$ of $L$ by $(W,\nabla)$ are equivalent
if there is an isomorphism $\phi:L_1\rightarrow L_2 $ of Lie-Rinehart algebras making the following diagrams commute:
\[
\diagram W \rto \dto^= & L_1 \rto \dto^{\phi} & L \dto^=    \\
         W \rto        & L_2 \rto             & L    
\enddiagram 
\]
Let $\Ext^1(L,W,\nabla)$ be the
set of equivalence classes of extensions of $L$ by the flat connection $(W,\nabla)$.
\end{definition}

Let $f\in \Z^2(L,W)$ be an element. It follows $f:L\times L\rightarrow W$ is $\B$-linear in both variables
with $f(x,x)=0$ for all $x\in L$ and $d^2(f)=0$.
Define the following product on $W\oplus L$:
\begin{align}
&\label{extension} [(w,x),(v,y)]:=(\nabla(x)(v)-\nabla(y)(w)+f(x,y), [x,y]) .
\end{align}
Let $L(f)$ be the left $\B$-module $W\oplus L$ equipped with the product $[,]$.
Define a map $\alpha_f:L(f)\rightarrow \Der_{\A}(\B)$ by $\alpha_f(w,x):=\alpha(x)$.
It follows the left $\B$-module $L(f)$ is a Lie-Rinehart algebra. The sequence
\[  W \rightarrow L(f) \rightarrow L  \]
 is an extension of $L$ by the flat connection $(W, \nabla)$.

Let $f,g\in \Z^2(L,W)$ be two 
cocycles. It follows there is an isomorphism $\phi: L(f)\rightarrow L(g)$ of extensions of Lie-Rinehart algebras 
if and only if there is a an element $\rho\in \C^1(L,W)$ with $d^1\rho = f-g$.
It follows we get a well defined map of sets
\[ \beta:\Z^2(L,W)\rightarrow \Ext^1(L,W,\nabla).\]
defined by sending $f$ to the equivalence class in $\Ext^1(L,W,\nabla)$ determined by $L(f)$.
Let $f+d^1\rho$ be an element in $\Z^2(L,W)$ with $\rho \in \C^1(L,W)$. It follows from the discussion above
that $\beta(f)=\beta(f+d^1\rho)$. We get a well defined map
\begin{align}
&\label{lie2} \overline{\beta}: \H^2(L,W) \rightarrow \Ext^1(L,W,\nabla) 
\end{align}
defined by
\[\overline{\beta}(\overline{f}):=L(f).\]

\begin{theorem} \label{h2} Let $\overline{\beta}$ be the map from (\ref{lie2}).
If $(L,\alpha)$ is an arbitrary Lie-Rinehart algebra the map $\overline{\beta}$ 
is an injection of sets. If $L$ is a projective $\B$-module it follows the map $\overline{\beta}$ is an 
isomorphism of sets.
\end{theorem}
\begin{proof} See  \cite{huebschmann}, Theorem 2.6.
\end{proof}

Note: One may construct an $\A$-module structure on $\Ext^1(L,W,\nabla)$ and one checks that the map $\overline{\beta}$ 
is an $\A$-linear map. Hence if $L$ is a projective $\B$-module there is an isomorphism
\[ \H^2(L,W)\cong \Ext^1(L,W,\nabla) \]
of left $\A$-modules.

One checks that
\[ \H^1(L,W)=\Der(L,W)/\Der^{inn}(L,W).\]

\begin{example} Cohomology of Lie algebras.\end{example}

The following result is well known from the cohomology theory of Lie algebras:

\begin{corollary} Let $L$ be a Lie algebra over  a field $k$ and let $W$ be a left $L$-module.
There is a bijection between $\H^2(L,W)$ and the set of equivalence classes of extensions of $L$
by $W$.
\end{corollary}
\begin{proof} The proof follows from Theorem \ref{h2} if we let $\A=\B=k$.
\end{proof}

\begin{example} \text{Singular cohomology of  complex algebraic manifolds.}

Assume $A$ is a finitely generated regular algebra over the complex numbers and let $X:=\Spec(A)$ be the associated
affine scheme. Let $X(\Cx)$ be the complex manifold associated to $X$ and let $L:=\Der_\Cx(A)$ be the Lie-Rinehart
algebra of derivations of $A$. It follows there is an isomorphism
\[ \H^i(L,A)\cong \H^i_{sing}(X(\Cx), \Cx) \]
of cohomology groups where $\H^i_{sing}(X(\Cx),\Cx)$ is singular cohomology of $X(\Cx)$ with complex coefficients.
It follows we get an isomorphism
\[ \Ext^1(L, A, \alpha) \cong \H^2_{sing}(X(\Cx),\Cx) \]
of complex vector spaces.
Hence to each cohomology class $\gamma \in \H^2_{sing}(X(\Cx),\Cx)$  we get an extension
\[ A \rightarrow L(\gamma) \rightarrow L \]
of Lie-Rinehart algebras. The class $\gamma$ is a purely topological object and the extension $L(\gamma)$ is a
purely algebraic object: $L(\gamma)$ is an infinite dimensional extension of the complex Lie algebra $L=\Der_\Cx(A)$
of $\Cx$-derivations of $A$.
\end{example}

\section{A PBW-Theorem for the twisted universal enveloping algebra}

In this section we generalize some constructions for Lie algebras and enveloping algebras of Lie algebras
from \cite{rinehart} and \cite{sridharan}  
to the case of Lie-Rinehart algebras and universal enveloping algebras of Lie-Rinehart algebras.
For an arbitrary Lie-Rinehart algebra $(L,\alpha)$ and an arbitrary cocycle $f\in \Z^2(L,\B)$, we define the 
universal enveloping algebra of type $f$ denoted $U(\B,L,f)$, and prove some basic properties of this algebra.
We prove a Poincare-Birkhoff-Witt Theorem for $U(\B,L,f)$ when $L$ is a projective $\B$-module, giving a simultaneous 
generalization of the Poincare-Birkhoff-Witt Theorem proved by Rinehart in \cite{rinehart} for Lie-Rinehart algebras
and Sridharan in \cite{sridharan} for Lie algebras (see Theorem \ref{pbw}).

Let for the rest of this section
$\alpha:L\rightarrow \Der_{\A}(\B)$ be a Lie-Rinehart algebra and let $f\in \Z^2(L,\B)$ be a cocycle.
Let $z$ be a generator for the free $\B$-module $F:=\B z$ and let
\[ F \rightarrow L(f) \rightarrow L \]
be the extension of $L$ by $F$ corresponding to $f$ as defined in (\ref{extension}).
For any elements $u:=az+x, v:=bz+y \in L(f)$ the following holds:
\[ [u,v]:=(\alpha(x)(b)-\alpha(y)(a)+f(x,y),[x,y]).\]
Write $x(b):=\alpha(x)(b)$.
The pair $(L(f),\alpha_f)$ where $\alpha_f(az+x):=\alpha(x)\in \Der_{\A}(\B)$ is by the results in the previous section
a Lie-Rinehart algebra. Hence $L(f)$ is a left $\B$-module and an $\A$-Lie algebra and $\alpha_f$ is a map
of left $\B$-modules and $\A$-Lie algebras.

Let $T(L(f)):=\oplus_{k\geq 0}L(f)^{\otimes_{\A} k}$ be the tensor algebra (over $\A$) of the $\A$-Lie algebra $L(f)$.
Let $T^r(L(f)):=\oplus_{k\geq r}L(f)^{\otimes_{\A} k}$ and let $T_r(L(f)):=\oplus_{k=0}^r L(f)^{\otimes_{\A} k}$.
Let $U_f$ be the two sided ideal in $T(L(f))$ generated by the set of elements
\[ u\otimes v -v\otimes u -[u,v] \]
with $u,v\in L(f)$. Let $U(L(f)):=T(L(f))/U_f$ be the universal enveloping algebra of the $\A$-Lie algebra $L(f)$.

Let $p:T(L(f))\rightarrow U(L(f))$ be the canonical map and let 
\begin{align}
&\label{uplus}U^+:=p(T^1(L(f))).
\end{align} 
Let 
\[ p_{\B}:\B\rightarrow U^+ \]
be defined by
\[ p_{\B}(b):=p(bz) \]
for all $b\in \B$.
Let 
\[ p_L:L\rightarrow U^+ \]
be defined by
\[ p_L(x):=p(x)\]
for $x\in L$
Let finally
\[ p_{L(f)}:L(f)\rightarrow U^+ \]
be defined by
\[ p_{L(f)}(w):=p(w) \]
for $w\in L(f)$.

\begin{definition}\label{classuniv}  Let $f\in \Z^2(L,\B)$ and let $U^+\subseteq U(L(f))$ be the algebra defined in
(\ref{uplus}). 
Let $J_f$ be the two sided ideal in $U^+$ generated by the following set:
\[ \{ p_{L(f)}(bw)-p_{\B}(b)p_{L(f)}(w) :\text{ where $b\in \B$ and $w\in L(f)$} \}. \]
Let $U(\B,L,f):=U^+/J_f$. By definition $U(\B,L,f)$ is an associative $\A$-algebra.
Let $U(\B,L,f)$ be the \emph{universal enveloping algebra of $(L,\alpha)$ of type $f$}.
Let $U(B,L):=U(B,L,0)$.
We say $U(B,L)$ is the \emph{universal enveloping algebra of $(L,\alpha)$}.
\end{definition}
The algebra $U(B,L)$ was first introduced by Rinehart in \cite{rinehart}.
The algebra $U(\B,L,f)$ is a simultaneous generalization of the universal enveloping algebra $U(\B,L)$ 
introduced by Rinehart and the twisted universal enveloping algebra $\lg_f$
of a Lie algebra $\lg$ introduced by Sridharan in \cite{sridharan}. If $f=0$ it follows $U(\B,L,0)=U(\B,L)$ and if $\B=\A$
and $\lg=L$ it follows $U(\B,L,f)=\lg_f$.

Let $p_1:T^1(L(f))\rightarrow U(\B,L,f)$ be the canonical map. Let $U^l(\B,L,f):=p(T^l(L(f)))$ and $U_l(\B,L,f):=p(T_l(L(f)))$.
We get a filtration
\[  \cdots \subseteq U^k(\B,L,f) \subseteq U^{k-1}(\B,L,f) \subseteq \cdots \subseteq U^1(\B,L,f)=U(\B,L,f) \]
called the \emph{descending filtration of} $U(\B,L,f)$.
We moreover get a filtration
\[ U_1(\B,L,f) \subseteq U_2(\B,L,f) \subseteq \cdots \subseteq U_k(\B,L,f) \subseteq \cdots \subseteq U(\B,L,f) \]
called the \emph{ascending filtration of} $U(\B,L,f)$.

Note: If $\rho \in \C^1(L,\B)$ is a cocycle it follows there is an isomorphism $L(f)\cong L(f+d^1\rho)$
of extensions. It follows there is an isomorphism
\[ U(\B,L,f)\cong U(\B,L,f+d^1\rho) \]
of filtered associative $\A$-algebras. We get for any cohomology class $ c\in \H^2(L,\B)$ a universal enveloping algebra
$U(\B,L,c):=U(\B,L,f)$ where $f$ is some element in $\Z^2(L,\B)$ representing the cohomology
class $c$. The $\A$-algebra $U(\B,L,c)$ is by the above discussion well defined up to isomorphism of filtered $\A$-algebras.
It moreover follows $U(\B,L)$ has a descending filtration $U^k(\B,L)$ and an ascending filtration $U_k(\B,L)$.

\begin{lemma} 
There is a one-to-one correspondence between the set of left $U(\B,L)$-modules and the set of flat $L$-connections.
\end{lemma}
\begin{proof} The Lemma follows from \cite{rinehart}.
\end{proof}

Let $\B w$ be the free rank one $\B$-module on the element $w$ and 
let $\tilde{L}:=\B w\oplus L(f)$ with the following Lie-product:
\[ [aw+u, bv+v]:= (u(b)-v(a))w+[u,v].\]
Here $u(b):=\alpha_f(u)(b)$ where $\alpha_f:L(f)\rightarrow \Der_{\A}(\B)$ is the anchor map of $L(f)$.
As left $\B$-module it follows $\tilde{L}=\B w\oplus \B z\oplus L$.
There is a canonical map
\[ \tilde{\alpha}:\tilde{L}\rightarrow \Der_{\A}(\B) \]
defined by
\[ \tilde{\alpha}(aw+bz+x):=\alpha(x)\]
and the pair $(\tilde{L},\tilde{\alpha})$ is a Lie-Rinehart algebra.
Let $U(\B, L(f))$ be the universal enveloping algebra of the pair $(L(f), \alpha_f)$ in the sense of Definition 
\ref{classuniv}.
Let
\[ r_1:T^1(\tilde{L})\rightarrow U(\B,L(f)) \]
be the canonical map.
We get a map
\[ r:\tilde{L}\rightarrow U(\B,L(f)) \]
defined by
\[ r(w):=r_1(w)\]
for $w\in \tilde{L}$. Let $z':=r(z)$ and $w':=r(w)$.
Let 
\begin{align}
&\label{univz} U(\B,L(f), z'):=U(\B,L(f))(z'-1) .
\end{align}
It follows $U(\B,L(f),z')$ has a descending filtration 
$U^k(\B,L(f),z')$ and an ascending filtration $U_k(\B,L(f),z')$.

\begin{theorem} \label{caniso} Let $L$ be a Lie-Rinehart algebra and let $f\in \Z^2(L,B)$.
Let $U(\B,L(f),z')$ be the algebra defined in (\ref{univz}).
There is a canonical  isomorphism of filtered $\A$-algebras and left $\B$-modules
\[ \phi: U(\B,L(f), z')\cong U(\B,L,f).  \]
\end{theorem}
\begin{proof} Define the map $\phi'$ as follows:
\[ \phi':T^1(\tilde{L})\rightarrow U(\B,L,f) \]
by
\[\phi'(aw+bz+x):=(a+b)z+x.\]
One checks $\phi'$ gives a well defined map 
\[ \phi: U(\B,L(f),z')\rightarrow U(\B,L,f) \]
of $\A$-algebras.
One shows $\phi$ has an inverse hence the first claim follows. The map $\phi$ maps the descending (resp. ascending) filtration
of $U(\B,L(f),z')$ to the descending (resp. ascending) filtration of $U(\B,L,f)$. The Theorem follows.
\end{proof}

Let $p_f:L(f)\rightarrow U(\B,L(f))$ be the canonical map of left $\B$-modules and let $U_k(A,L(f))$ be the ascending
filtration of $U(A,L(f))$.

\begin{proposition} \label{generated} The abelian group $U_k(\B,L(f))$ is generated as left $\B$-module by the set
\[ \{ p_f(x_{i_1})p_f(x_{i_2})\cdots p_f(x_{i_l}) :\text{ with $x_{i_j}\in L(f)$ and $l\leq k$.}  \} \]
\end{proposition}
\begin{proof} We prove the result by induction in $k$. For $k=1$ it is obvious. Assume the result is true for the case $i=k-1$.
Assume $i=k$. Let $p=p_f$ and let $w=p(z_1)\cdots p(z_k)\in U_k(\B,L(f))$ with $z_i\in L(f)$. We get by the induction
hypothesis the following equality:
\[ p(z_2)\cdots p(z_k)=\sum_I a_I p(x_{i_1})\cdots p(x_{i_l}) \]
with $a_I\in \B$ and $x_{i_j}\in L(f)$ for all $I, i_j$. We may write $z_1=az+x\in L(f)$. We get
\[ p(z_1)p(z_2)\cdots p(z_k)= \sum_I (az+x)a_I p(x_{i_1})\cdots p(x_{i_l}) =\]
\[ \sum_I aa_Ip(x_{i_1})\cdots p(x_{i_l})+a_Ip(x)p(x_{i_1})\cdots p(x_{i_l}) +\alpha(x)(a_I)p(x_{i_1})\cdots p(x_{i_l}) \]
hence the claim holds for $i=k$. The Lemma follows.
\end{proof}

\begin{corollary} There is a canonical surjective map of left $\B$-modules
\[ \phi: \Sym^k_{\B}(L(f))\rightarrow U_k(\B,L(f))/U_{k-1}(\B,L(f)).\]
\end{corollary}
\begin{proof}  Assume $x_1,\ldots ,x_k\in L(f)$. By induction one proves the following result: Assume $\sigma$ is a permutation 
of the set $\{1,2,..,k\}$. The following formula holds:
\[ p(x_1)\cdots p(x_k)=p(x_{\sigma(1)})\cdots p(x_{\sigma(k)}) +w \]
with $w\in U_{k-1}(\B,L(f))$.
Define the following map:
\[\phi: \Sym^k_{\B}(L(f))\rightarrow U_k(\B,L(f))/U_{k-1}(\B,L(f)) \]
by
\[ \phi(x_1\cdots x_k):=\overline{p(x_1)\cdots p(x_k)}.\]
It follows
\[ \phi(x_1\cdots x_k)=\phi(x_{\sigma(1)}\cdots x_{\sigma(k)}) \]
hence $\phi$ is well defined. By Proposition \ref{generated} it follows the map $\phi$ is a surjective map of left $\B$-modules
and the Corollary is proved.
\end{proof}

\begin{lemma}\label{sym} Assume $L(f)$ is a projective $\B$-module. For all $k\geq 1$ there is a canonical isomorphism
of left $\B$-modules
\[ \Sym^k_{\B}(L) \cong         U_k(\B,L(f),z')/U_{k-1}(\B,L(f), z') .\]
\end{lemma}
\begin{proof} Let $p_f:L(f)\rightarrow U(\B,L(f))$ be the canonical map and let $z':=p_f(z)$. Recall that
$L(f):=\B z\oplus L$ where $z$ is a generator for the free rank one submodule $\B z$ of $L(f)$.
The element $z'$ is a central element in $U(\B,L(f))$: For all elements $w\in U(\B,L(f))$ it follows that
$z'w=wz'$. It follows $(z'-1)w=w(z'-1)$ for all $w\in U(\B,L(f))$. It follows the two sided ideal in $U(\B,L(f))$ generated
by $z'-1$ is the following set:
\[ \{ w(z'-1) :\text{ where $w\in U(\B,L(f))$.} \} .\]
We get a commutative diagram of exact sequences of left $\B$-modules
\[
\diagram 0 \rto & U_k(\B,L(f))(z'-1) \rto & U_k(\B,L(f)) \rto & U_k(\B,L(f),z') \rto & 0 \\
 0 \rto & U_{k-1}(\B,L(f))(z'-1) \rto \uto^u & U_{k-1}(\B,L(f)) \rto \uto^v & U_{k-1}(\B,L(f),z') \rto \uto^w & 0 
\enddiagram .
\]
Since $ker(u)=ker(v)=ker(w)=0$ we get by the snake lemma a short exact sequence of left $\B$-modules
\[ 0 \rightarrow coker(u) \rightarrow^i coker(v) \rightarrow^j coker(w) \rightarrow 0.\]
There is by definition an isomorphism of left $\B$-modules
\[ coker(w)\cong U_k(\B,L(f),z')/U_{k-1}(\B,L(f),z').\]
By assumption there is a canonical isomorphism of left $\B$-modules
\[  \Sym^k_{\B}(L(f))\cong U_k(\B,L(f))/U_{k-1}(\B,L(f)) .\]
There is also an isomorphism
\[ \Sym^k_{\B}(L(f))\cong \Sym^{k-1}_{\B}(L(f))z\oplus \Sym^k_{\B}(L).\]
One checks that $im(i)=\Sym^{k-1}_{\B}(L(f))z$ hence we get an isomorphism
\[ \Sym^k_{\B}(L)\cong coker(w)\cong U_k(\B,L(f),z')/U_{k-1}(\B,L(f),z'),\]
and the Lemma is proved.
\end{proof}

\begin{theorem}\label{pbw}(The PBW-Theorem) Let $\alpha:L\rightarrow \Der_R(A)$ be a Lie-Rinehart algebra 
where $L$ is a projective $\B$-module. There is a canonical isomorphism of graded $\B$-algebras
\[ \Sym^*_{\B}(L)\cong Gr(U(\B,L,f)), \]
where $Gr(U(\B,L,f))$ is the associated graded algebra of $U(\B,L,f)$ with respect to the ascending filtration.
\end{theorem}
\begin{proof} The Theorem follows from Theorem \ref{caniso} and Lemma \ref{sym}
\end{proof}

Note: When $f=0$ is the zero cocycle we get the following result: There is a canonical isomorphism
of graded $\B$-algebras
\[ \Sym^*_{\B}(L)\cong Gr(U(\B,L))\]
where $U(\B,L)$ is Rineharts enveloping algebra of the Lie-Rinehart algebra $L$.
When $\A=\B$ and $\lg=L$ we get the following result:
There is a canonical isomorphism of graded $\B$-algebras
\[ \Sym^*_{\B}(\lg)\cong Gr(\lg_f) \]
where $\lg_f$ is Sridharan's twisted universal enveloping algebra of the $\B$-Lie algebra $\lg$.
Hence Theorem \ref{pbw} gives a simultaneous generalization of results of 
Rinehart and Sridharan obtained in \cite{rinehart} and \cite{sridharan}.

\begin{corollary} Assume $\B$ is a Noetherian ring and $L$ a finitely generated projective $\B$-module.
Let $f\in \Z^2(L,\B)$. It follows $U(\B,L,f)$ is Noetherian.
\end{corollary}
\begin{proof} By Theorem \ref{pbw} it follows $Gr(U(\B,L,f)\cong \Sym^*_{\B}(L)$ and by hypothesis
$\Sym^*_{\B}(L)$ is Noetherian. It follows $U(\B,L,f)$ is Noetherian.
\end{proof}

\begin{example} \text{The classical PBW-Theorem.}

Let $\lg$ be a finite dimensional Lie algebra over a field $k$ with basis $B=\{e_1,..,e_k\}$ and let
$B_U:=\{ e_1^{p_1}\cdots e_k^{p_k}: p_i\geq 0.\}$. Let $U(\lg)$ be the universal enveloping algebra of the $k$-Lie algebra $\lg$
and view $B_U$ as a subset of $U(\lg)$.
There is a canonical map

\begin{align}
&\label{pbwlie} \gamma: \Sym^*_k(\lg) \rightarrow  Gr(U(\lg)) 
\end{align}

of graded $k$-algebras. If the set $B_U$ generates $U(\lg)$ as left $k$-module it follows the map $\gamma$ is surjective.
If the set $B_U$ is linearly independent over the field $k$ it follows the map $\gamma$ is injective. Hence
the classical PBW-theorem for $U(\lg)$ is equivalent to the fact that $\gamma$ is an isomorphism of graded $k$-algebras.

The PBW-Theorem for Lie-Rinehart algebras is a globalized version of (\ref{pbwlie})
valid for a Lie-Rinehart algebra $L$ which is projective as
$B$-module. 

Assume $X:=\Spec(\B)$ is an affine scheme with structure sheaf $\O$. Let $h:\A\rightarrow \B$ 
be a ring homomorphism and let $S:=\Spec(\A)$ and $\pi:X\rightarrow S$ the canonical map.
Let $\alpha:L\rightarrow \Der_{\A}(\B)$ be a Lie-Rinehart algebra which is finitely generated and projective as left $\B$-module, 
and let $\mathcal{L}$ be the $\O$-module associated to $L$. Let $\mathcal{U}_{L,f}$ be the $\O$-module associated to $U(\B,L,f)$.
We view $\mathcal{U}_{L,f}$ as a sheaf of filtered associative unital $\pi^{-1}(\O_S)$-algebras on $X$.
The canonical isomorphism
\[ \gamma: \Sym_{\B}^*(L) \cong Gr(U(\B,L,f)) \]
implies the following: Around every point $x\in X$ there is a Zariski open subset $U$
where the sections of $\mathcal{L}$ over $U$ has a generating
set $s_1,..,s_k$ as free $\O(U)$-module with the following property: Let
$B_U:=\{s_1^{p_1}\cdots s_k^{p_k}:p_i\geq 0\}$ and view $B_U$ as a subset of the left $\O(U)$-module $\mathcal{U}_{L,f}(U)$.
Then the fact that $\gamma$ is an isomorphism implies that
the set $B_U$ is a linearly independent set of generators of $\mathcal{U}_{L,f}(U)$
viewed as left $\O(U)$-module. This property holds for any point $x\in X$. Hence we may view $\mathcal{U}_{L,f}$ as an infinite
rank locally free sheaf on $X$.
\end{example}

\section{Application I: Deformations of almost commutative rings}

In this section we study the deformation groupoid $\Adef (\Sym_{\B}^*(L))$ of a Lie-Rinehart algebra $L$ which is 
projective as left $\B$-module. We show that isomorphism classes of objects 
in $\Adef (\Sym_{\B}^*(L))$ are parametrized by the cohomology
group $\H^2(L,\B)$ and that the morphisms in $\Adef (\Sym_{\B}^*(L))$ are parametrized by the group $\Z^1(L,\B)$
(see Theorem \ref{maindeform} and \ref{morphisms}). An application of these results is that for any filtered almost
commutativ ring $U$ satisfying a PBW-condition it follows that the category of left $U$-modules is equivalent
to the category of $L$-connections of curvature type $f$, where $L$ is a Lie-Rinehart algebra and $f\in \Z^2(L,\B)$ is
a 2-cocycle (see Corollary \ref{modules}).

The deformation groupoid $\Adef (\Sym_{\B}^*(L))$ was first introduced by Sridharan in the case of a Lie-Algebra $\lg$ 
where $\lg$ is a free module over a commutative ring $A$ (see \cite{sridharan}).

\begin{definition}
Let $(U,U_i)$ be a filtered associative algebra with filtration
\[ U_0\subseteq U_1\subseteq \cdots \subseteq U_k \subseteq \cdots U \]
where $U_0=\B$ and $h:\A\rightarrow \B$ is an arbitrary injective map of commutative rings with unit. 
Assume $\A\subseteq Center(U)$
and assume moreover that for any elements $x\in U_i, y\in U_j$ it follows
$xy\in U_{i+j}$.
We say $(U,U_i)$ is \emph{almost commutative} if the following holds: Assume 
$x_1\in U_{i_1},..,x_k\in U_{i_k}$ and  assume $\sigma$ is a permutation of $k$ elements. Let $i=:\sum_{l=1}^k i_l$.
 There is an equality
\[ x_1\cdots x_k =x_{\sigma(1)}\cdots x_{\sigma(k)}+y_{i-1} \]
where $y_{i-1}\in U_{i-1}$. 
\end{definition}

It follows the associated graded algebra $Gr(U)$ of $(U,U_i)$ is commutative.

\begin{example} \text{Rings of differential operators.}
Let $\D_{\A}^0(\B):=\B$ and let
\[ \D_{\A}^k(B):=\{ \partial \in \End_{\A}(\B): [\partial , a]\in \D_{\A}^{k-1}(\B) \text{ for all $a\in \B$} \}. \]
Let $\D_{\A}(\B):=\cup_{k\geq 0}\D_{\A}^k(\B)$ be the \emph{ring of differential operators of $\B/\A$}. 
The associative ring  $\D_{\A}(\B)\subseteq \End_{\A}(\B)$ has a filtration
\begin{align}
&\label{diff_filt}\B:=\D_{\A}^0(\B)\subseteq \D_{\A}^1(\B)\subseteq \cdots \subseteq 
\D_{\A}^k(\B)\subseteq \cdots \subseteq \D_{\A}(\B) .
\end{align}
The filtered ring $(\D_{\A}(\B), \D_{\A}^k(\B) )$ is almost commutative.
\end{example}

Assume $(U,U_i)$ is almost commutative and let  $L:= U_1/U_0$.
Let
\begin{align}
&\label{gammaU} \gamma_U:\Sym^*_{\B}(L)\rightarrow Gr(U) 
\end{align}
be the canonical map of graded $\B$-algebras.

\begin{definition} Let $(U,U_i)$ be an almost commutative algebra. We say $(U,U_i)$ is  \emph{almost commutative of PBW-type} if 
the canonical map $\gamma_U$ from (\ref{gammaU}) is an isomorphism of graded $\B$-algebras.
\end{definition}

\begin{example} \label{diffII}\text{Rings of differential operators II.}

Let $L:=\D^1_{\A}(\B)/\D^0_{\A}(\B)$ and consider the canonical map
\[ \gamma:\Sym^*_{\B}(L)\rightarrow Gr(\D_{\A}(\B)).\]
The map $\gamma$ is neither surjective nor injective in general. If $\A=k$ is a field of characteristic zero
and $\B$ is a regular $k$-algebra of finite type it follows the map $\gamma$ is surjective. This is because
the algebra $\D_k(\B)$ is generated by $\Der_k(\B)$: Every polynomial 
differential operator $\partial:\B\rightarrow \B$ may be written as a sum of products of derivations.
\end{example}

Assume in the following that $(U,U_i)$ is almost commutative of PBW-type.

We get an exact sequence of left $\B$-modules
\[  U_0 \rightarrow U_1 \rightarrow L  .\]
Consider the following map
\[\psi: U_0\times U_1 \rightarrow L\]
where
\[ \psi(b,z):=b\overline{z}\]
where $\overline{z}\in L=U_1/U_0$ is the equivalence class of $z$. Since $U$ is an associative algebra
it follows $U_1$ is a left and right $\B$-module and since $\Sym_{\B}^*(L)$ is a commutative $\B$-algebra it follows
the element $b\overline{z}-\overline{z}b$ is zero in $L$. It follows the commutator $[z,b]=zb-bz$ is an element in
$U_0\subseteq U_1$.
We get a map
\[ \tilde{\alpha}:U_1\rightarrow \End(B) \]
defined by
\[ \tilde{\alpha}(z)(b):=[z,b].\]
It follows immediately that $\tilde{\alpha}(z)\in \End_{\A}(\B)$
for any element $z\in U_1$. We moreover get the following equation:
\[ \tilde{\alpha}(z)(ab)=[z,ab]=zab-azb+azb-abz=[z,a]b-a[z,b]=\tilde{\alpha}(z)(a)b+a\tilde{\alpha}(z)(b)\]
hence
\[\tilde{\alpha}(z)\in \Der_{\A}(\B).\]
It follows we get a map
\begin{align}
&\label{alphat} \tilde{\alpha}:U_1\rightarrow \Der_{\A}(\B).
\end{align}
\begin{lemma} Assume $(U,U_i)$ is almost commutative of PBW-type and assume $\tilde{\alpha}$ is the 
map from (\ref{alphat}). The pair $\{U_1,\tilde{\alpha}\}$ is a Lie-Rinehart algebra.
\end{lemma}
\begin{proof} The proof is an exercise.
\end{proof}

Since $U_0\subseteq U_1$ is an ideal we get an induced structure of $\A$-Lie algebra on $L=U_1/U_0$.
By definition $\B=U_0\subseteq U_1$ is an abelian sub-algebra. It follows the natural sequence
\[ \B \rightarrow U_1 \rightarrow L \]
is an exact sequence of Lie-Rinehart algebras. We get an induced Lie-Rinehart structure
\begin{align}
&\label{alpha}  \alpha:L\rightarrow \Der_{\A}(\B) .
\end{align}
\begin{definition} \label{filteredalgebra} Assume $(U,U_i)$ is almost commutative of PBW-type
with induced anchor map $\alpha:L\rightarrow \Der_{\A}(\B)$ from (\ref{alpha}), where $L=U_1/U_0$.
We say $(U,U_i)$ is a \emph{filtered algebra of type $\alpha$}.
Let $c(U)\in \Ext^1(L,\B,\alpha)$ be the \emph{ deformation class} defined by the extension
\[  \B\rightarrow U_1 \rightarrow L  .\]
We say $U$ is the \emph{trivial deformation} if $c(U)=0$ in $\Ext^1(L,\B,\alpha)$
\end{definition}

\begin{example}\text{Rings of differential operators III.}
The assumptions are as in Example \ref{diffII}. The canonical map
\[ \gamma: \Sym_A(L)\rightarrow Gr(\D_R(A)) \]
is seldom an isomorphism, hence $(\D_R(A), \D^k_R(A))$ is seldom a filtered algebra of type
$\alpha$ for some Lie-Rinehart algebra $\alpha: L\rightarrow \Der_R(A)$.

\end{example}

\begin{example} \text{Twisted rings of differential operators.}

If $A$ is a regular algebra of finite type over the complex numbers, 
$L=\Der_{\Cx}(A)$ and $f\in \Z^2(L,A)$ it follows $U(A,L,f)$ is an algebra of twisted differential operators
on $A$ in the sense of \cite{beilinson}, Definition 2.1.1.

\end{example}

\begin{example} \text{Simpson's sheaf of rings of differential operators.}

Recall the following axioms from Simpson's paper (see \cite{simpson}): Let $X$ be a smooth algebraic variety of finite type
over the complex numbers and let $\Gamma$ be a sheaf of filtered associative unital algebras on $X$ with filtration
$\Gamma_i$, where $i=0,1,2..$. Here $\Gamma_i\subseteq \Gamma$ is a sub-sheaf of abelian groups. Let $K=\mathbb{C}_X$ 
be the constant sheaf on $X$ on the field $\mathbb{C}$ of complex numbers.
Simpson defines $\Gamma$ to be a \emph{sheaf of rings of differential operators} on $X$ relative to $K$ if the following holds:
\begin{align}
&\label{si1} \text{The sheaf $K$ is in the center of $\Gamma$.}\\
&\label{si2} \Gamma_0=\O_X, \Gamma_i\Gamma_j\subseteq \Gamma_{i+j}.\\
&\label{si3}\text{The left and right $\O_X$-module structure on $\Gamma_i/\Gamma_{i-1}$ coincide.}\\
&\label{si4}\text{The graded $\O_X$-modules $\Gamma_i/\Gamma_{i-1}$ are coherent.}\\
&\label{si5}\text{The canonical map $\gamma:\Sym_{\O_X}(\Gamma_1/\Gamma_0)\rightarrow Gr(\Gamma,\Gamma_i)$ is an isomorphism.}
\end{align}

If $X=\Spec(\B)$ is a finitely generated regular algebra over the complex numbers one checks that the sheaf of filtered algebras
$(\Gamma, \Gamma_i)$ satisfying \ref{si1} - \ref{si5} corresponds to a filtered algebra $(U,U_i)$ of type $\alpha$ 
from Definition \ref{filteredalgebra} if we let $U=\H^0(X, \Gamma)$ and $U_i=\H^0(X, \Gamma_i)$.
\end{example}

\begin{example} \text{The Riemann-Hilbert correspondence.}

In Simpson's paper \cite{simpson} the Riemann-Hilbert correspondence is studied at the level of moduli spaces. The 
Riemann-Hilbert correspondence
is the following: Let $X$ be a smooth projective variety of finite type over the complex numbers. There is an equiuvalence 
of categories between the category of finite rank complex representations of the topological fundamental group of $X$ and the 
category of finite rank complex vector bundles on $X$ with a flat connection. There is moreover an equivalence of 
categories between the category of finite rank complex representations of the fundamental group of $X$ and the category of 
Higgs-bundles on $X$. This gives homeomorphisms between
the moduli space of representations of the topological fundamental group of $X$, 
the moduli space of finite rank complex vector bundle  on $X$ with a flat connection and 
the moduli space of Higgs-bundles on $X$. The fact that $X$ is and algebraic variety 
induce algebraic structures on these moduli spaces and Simpson invesitgates the relationship between these algebraic 
structures and the underlying topological equivalences. In this paper all techniques are algebraic, 
giving a theory valid in the affine case over any base ring. 
\end{example}

Assume now $(U,U_i)$ is a filtered algebra of type $\alpha$ where $L=U_1/U_0$ and $\alpha:L\rightarrow \Der_{\A}(\B)$.
Assume moreover that $L$ is a projective $\B$-module.  Consider the exact sequence
\begin{align}
&\label{elie} U_0 \xrightarrow{i} U_1 \xrightarrow{p}  U_1/U_0  
\end{align}
of Lie-Rinehart algebras.
Assume $t$ is a right splitting of \ref{elie}, hence $t:U_1/U_0 \rightarrow U_1$ is left $B$-linear and $p\circ t=id$.
Let
\[ \phi_{U,1}:L\rightarrow U_1/U_0 \]
be the first component of the graded isomorphism $\phi_U:\Sym_{\B}^*(L)\cong Gr(U)$.
Let $\phi_{U,1}^{-1}$ be the inverse and let $T=t\circ \phi_{U,1}$ and $P=p\circ \phi_{U,1}^{-1}$.
We get an exact sequence
\[ U_0 \rightarrow U_1 \rightarrow^P L  \]
which is right split by $T$.

Assume $p(z)=x$ and let $\gamma:L\rightarrow \Der_{\A}(\B)$ be defined by
\[\gamma(x)(b):=[T(x), b]=T(x)b-bT(x).\]
Assume $(U,U_i)$ is a filtered algebra of type $\alpha$. This means that
\[ \gamma(x)(b):=[T(x),b]=T(x)b-bT(x)=\alpha(x)(b).\]
Assume moreover that
\[ [T(x),T(y)]-T([x,y])=f(x,y) \in \B \subseteq U_1\]
where $f\in \Z^2(L,\B)$.

\begin{lemma} There is a commutative diagram
\[
\diagram   Gr(U(\B,L,f)) \rrto^{Gr(\tilde{T})} &  &  Gr(U) \\
           \Sym_{\B}^*(L)\uto^{\phi_f} \urrto^{\phi_U} & \\
\enddiagram.
\]
\end{lemma}
\begin{proof} 
Recall the construction of the algebra $U(\B,L,f)$. Let $L(f):=\B z\oplus L$ with product defined by (\ref{extension}).
Recall the canonical map
\[\sigma_1:T^1(L(f))\rightarrow U(\B,L,f).\]
Define
\[T':T^1(L(f))\rightarrow U\]
by
\[ T'((a_1z+x_1)\otimes \cdots \otimes (a_kz+x_k)):=\prod_i(a_i+T(x_i)).\]
It follows
\[ T'((az+x)\otimes (bz+y)-(bz+y)\otimes (az+x)-[az+x,bz+y])=\]
\[ (a+T(x))(b+T(y))-(b+T(y))(a+T(x))- (\alpha(x)(b)-\alpha(y)(a)+f(x,y))z-T([x,y])=\]
\[ab+aT(y)+T(x)b+T(x)T(y)-ba-bT(x)-T(y)a-T(y)T(x)-\alpha(x)(b)+\]
\[\alpha(y)(b)-f(x,y)-T([x,y])=0\]
since $T(x)b-bT(x)=\alpha(x)(b)$.
Moreover for any $b\in \B$ and $w:=az+x \in L(f)$ it follows
\[ T'(\sigma_1(bw)-\sigma_1(b)\sigma_1(w))=T'(baz+bx-bzaz-bzx)=0\]
hence $T'$ induce a map
\[ \tilde{T}:U(\B,L,f)\rightarrow U\]
of filtered algebras:
\[ \tilde{T}(x_1\cdots x_k)=T(x_1)\cdots T(x_k)=\overline{t(\phi_{U,1}(x_1))\cdots t(\phi_{U,1}(x_k))} \]
for $x_i\in L$.
Since $p\circ t \circ \phi_{U,1}=\phi_{U,1}=\overline{t\circ \phi_{U,1}}$ it follows
\[ \tilde{T}(x_1\cdots x_k)=\phi_{U,1}(x_1)\cdots \phi_{U,1}(x_k).\]
The Lemma follows.
\end{proof}

Hence there is an equality $Gr(\tilde{T})\circ \phi_f = \phi_U$ hence $Gr(\tilde{T})=\phi_U\circ \phi_f^{-1}$.
It follows the map
\[ Gr(\tilde{T}):Gr(U(\B,L,f))\rightarrow Gr(U)\]
is an isomorphism of filtered algebras.
\begin{lemma} \label{isorings} The map $\tilde{T}:U(B,L,f)\rightarrow U$ is an isomorphism of associative rings.
\end{lemma}
\begin{proof} Since $Gr(\tilde{T})$ is an isomorphism it follows the induced map
\[ \tilde{T}:U_0(\B,L,f)\rightarrow U_0 \]
is an isomorphism. Assume the induced map
\[ \tilde{T}:U_{k-1}(\B,L,f)\rightarrow U_{k-1} \]
is an isomorphism. We get a commutative diagram of exact sequences
\[
\diagram  0   \rto  &   U_{k-1}(\B,L,f) \rto \dto^{\tilde{T}}   &    U_k(\B,L,f) \rto \dto^{\tilde{T}} &
U_k(\B,L,f)/U_{k-1}(\B,L,f) \rto  \dto^{Gr(\tilde{T})_k} & 0 \\
0 \rto   &    U_{k-1}  \rto & U_k \rto & U_k/U_{k-1} \rto & 0
\enddiagram.
\]
It follows from the snake Lemma that the induced morphism
\[ \tilde{T}: U_k(\B,L,f)\rightarrow U_k \]
is an isomorphism. The Lemma follows by induction.
\end{proof}

\begin{definition}\label{deformationgroupoid}
Let $\Adef (\Sym_{\B}^*(L))$ be the following category: Let the objects of $\Adef (\Sym_{\B}^*(L))$ be the set of pairs
$(U,\psi_U)$, where $U$ is a filtered algebra of type $\alpha$ and where
\[ \psi_U: \Sym_{\B}^*(L)\rightarrow Gr(U) \]
is a fixed isomorphism of graded $\B$-algebras.
A morphism $\theta: (U,\psi_U)\rightarrow (V,\psi_V)$ in $\Adef (\Sym_{\B}^*(L))$ is a map
of filtered algebras
\[ \theta: U\rightarrow V\]
such that the induced map on associated graded rings
\[ Gr(\theta):Gr(U)\rightarrow Gr(V) \]
satisfies $Gr(\theta)\circ \psi_U=\psi_V$. Since $\psi_U$ and $\psi_V$ are isomorphisms, it follows 
\[ Gr(\theta)=\psi_V\circ \psi_U^{-1} ,\]
hence the map $Gr(\theta)$ is an isomorphism of graded $\B$-algebras. It follows the map $\theta$
is an isomorphism of filtered algebras. One checks the inverse $\theta^{-1}$ is a map in $\Adef (\Sym_{\B}^*(L))$, hence
the category $\Adef (\Sym_{\B}^*(L))$ is a groupoid. 
The category $\Adef (\Sym_{\B}^*(L))$ is the \emph{deformation groupoid of $(L,\alpha)$}.
\end{definition}

The category $\Adef (\Sym_{\B}^*(L))$ was first introduced in \cite{sridharan} for Lie algebras
over rings.

Let $\Iso (\Adef (\Sym_{\B}^*(L)))$ be the set of isomorphism classes of 
objects in $\Adef (\Sym_{\B}^*(L))$ and define the following  map:
\begin{align}
&\label{groupoid_equivalence} h: \H^2(L,\B) \rightarrow \Iso (\Adef (\Sym_{\B}^*(L))) 
\end{align}
by
\[ h(\tilde{f}):=( U(\B,L,f), \phi_f). \]
Here
\[ \phi_f:\Sym_{\B}^*(L)\rightarrow Gr(U(\B,L,f)) \]

is the canonical isomorphism of graded $\B$-algebras defined above. 
The map is well defined since for two elements $f,f+d^1\rho$  representing
the cohomology class $\tilde{f}$ in $\H^2(L,\B)$ it follows there is an isomorphism
\[ U(\B,L,f)\cong U(\B,L,f+d^1\rho)\]
of filtered algebras in $\Adef (\Sym_{\B}^*(L))$.

\begin{theorem}\label{maindeform}
  Assume $\alpha:L\rightarrow \Der_{\A}(\B)$ is a Lie-Rinehart algebra which is projective as left $\B$-module.
Let the map
\[ h: \H^2(L,\B) \rightarrow \Iso (\Adef (\Sym_{\B}^*(L))) \]
be the map defined in (\ref{groupoid_equivalence}).
It follows the map $h$ is a one to one correspondence of sets.
\end{theorem}
\begin{proof} By Lemma \ref{isorings} it follows $h$ is a surjective map. Assume $h(f)=h(g)$ for two elements
$f,g\in \Z^2(L,\B)$. It follows we get an isomorphism
\[ U(\B,L,f)\cong U(\B,L,g) \]
of filtered algebras.
It follows we get isomorphic extensions of Lie-Rinehart algebras $L(f)\cong L(g)$, hence there is an element
$\rho \in C^1(L,\B)$ with $d^1\rho=f-g$. Hence $\tilde{f}=\tilde{g}$ in $\H^2(L,\B)$ and the map is injective. 
The Theorem is proved.
\end{proof}

Theorem \ref{maindeform} was first proved in \cite{sridharan} for Lie algebras over an arbitrary base ring $K$.
The paper \cite{sridharan} use the classical Chevalley-Eilenberg complex to define the cohomology
of a Lie algebra. The Lie-Rinehart complex of a Lie-Rinehart algebra $\alpha:L\rightarrow \Der_{\A}(\B)$ specialize
to the Chevalley-Eilenberg complex if we let $\A=\B$.  
In \cite{tortella} a similar result is proved for a holomorphic Lie algebroid on a
compact K\"{a}hler-manifold $X$ using complex analytic techniques. If $X$ is projective it follows $X$ is algebraic.
Any holomorphic finite rank vector bundle on $X$ is algebraic and a holomorphic Lie algebroid corresponds to
a sheaf of Lie-Rinehart algebras. Hence in the complex projective situation we may use algebraic techniques 
and a global version of Lie-Rinehart cohomology to get a similar classification.
The cohomology theory used in \cite{tortella} is defined using a global complex analytic version of 
the Lie-Rinehart complex. The result is also mentioned (without proof) in \cite{beilinson} in Section 2.1.13.

There is work in progress on a globalization of the results in this section to give results for an arbitrary scheme
$X$ over a fixed base scheme $S$.

\begin{example} \text{A classification of sheaves of rings of differential operators.}

Let $\pi:\Spec(\B)\rightarrow \Spec(\A)$ be an arbitrary morphism of affine schemes and let $X:=\Spec(\B), S:=\Spec(\A)$. Let
$(\Gamma, \Gamma_i)$ be a filtered sheaf of associative rings satisfying axiom \ref{si1} - \ref{si5} with the property
that $\H^0(X, \Gamma_1/\Gamma_0)$ is a projective $\H^0(X,\Gamma_0)$-module. 
It follows from Theorem \ref{maindeform} there is a Lie-Rinehart algebra 
\[ \alpha:L\rightarrow \Der_{\A}(\B) \]
together with an element $f\in \Z^2(L,\B)$ and an isomorphism $\H^0(X, \Gamma)\cong U(\B,L,f)$ of filtered associative rings.
Hence the sheafification of $U(\B,L,f)$ gives rise to a sheaf $\mathcal{U}_{L,f}$ 
of differential operators on $X$ relative to $S$ isomorphic to $\Gamma$. When we vary $L$ and $f$ we get all such sheaves.
\end{example}

Let $\mathbb{V}(\H^2(L,\B)):=\Spec(\Sym_{\A}^*(\H^2(L,\B)^*))$. It follows
\[\pi:\mathbb{V}(\H^2(L,\B))\rightarrow \Spec(\A) \]
is a scheme over $\A$. If $\H^2(L,\B)$ is a locally free $\A$-module it follows $\mathbb{V}(\H^2(L,\B))$ is a vector bundle 
over $\Spec(\A)$. Theorem \ref{maindeform} shows that isoclasses in $\Adef (\Sym_{\B}^*(L))$ are parametrized by the points 
of the scheme $\mathbb{V}(\H^2(L,\B))$. 

\begin{definition}
Let $\nabla:L\rightarrow \End_{\A}(W)$ be an $L$-connection and let $f\in \Z^2(L,\B)$.
We say $\nabla$ is an \emph{L-connection of curvature type $f$} if the following is satisfied:
For all $x,y\in L$ and $v\in W$ the following formula holds:
\[ R_\nabla(x\wedge y)(v)=f(x,y)v.\]
Here $R_\nabla$ is the curvature of $\nabla$.
\end{definition}

\begin{lemma} \label{conncurvf} Let $W$ be a left $\B$-module.
There is a one-to-one correspondence
between the set of $L$-connections of curvature type $f$ on $W$ and the set of flat $L(f)$-connections
on $W$ with $\nabla(z)=Id_W$.
\end{lemma}
\begin{proof} Given an $L$-connection $\nabla:L\rightarrow \End_{\A}(W)$ of curvature type $f$ with $f\in \Z^2(L,\B)$.
Define the following map:
\[ \onb:L(f)\rightarrow \End_{\A}(W)\]
by
\[ \onb(az+x):=aId_W+\nabla(x).\]
One checks the following holds: $\onb(u+v)=\onb(u)+\onb(v)$ and $\onb(au)=a\onb(u)$ for all $u,v\in L(f)$ and $a\in \B$.
Hence $\onb$ is an $\B$-linear map. Moreover one checks that for any $a\in \B, u\in L(f)$ and $w\in W$ the following holds:
\[ \onb(u)(aw)=a\onb(u)(w)+\alpha_f(u)(a)w.\]
Hence $\onb$ is an $L(f)$-connection on $W$.
We calculate the curvature of $\onb$. Let $u:=(a,x), v:=(b,y)\in L(f)$. By definition it follows
\[ R_{\onb}(u,v):=[\onb(u),\onb(v)]-\onb([u,v])=\]
\[ [aId_W+\nabla(x),bId_W+\nabla(y)]-\onb((\alpha(x)(b)-\alpha(y)(a)+f(x,y),[x,y])).\]
We get
\[[aId_W+\nabla(x),bId_W+\nabla(y)]=\]
\[\nabla(x),\nabla(y)]-\alpha(y)(a)Id_W+\alpha(x)(b)Id_W.\]
We also get
\[ \onb((\alpha(x)(b)-\alpha(y)(a)+f(x,y),[x,y]))=\]
\[(\alpha(x)(b)-\alpha(y)(a)+f(x,y))Id_W+\nabla([x,y]).\]
It follows
\[ R_{\onb}(u,v)=\]
\[[\nabla(x),\nabla(y)]-\alpha(y)(a)Id_W+\alpha(x)(b)Id_W- ((\alpha(x)(b)-\alpha(y)(a)+f(x,y))Id_W+\nabla([x,y]))=\]
\[R_{\nabla}(x,y)-f(x,y)Id_W=0\]
since $\nabla$ has curvature type $f$. 
It follows $\onab$ is a  flat $L(f)$-connection.

Assume $\onb:L(f)\rightarrow \End_{\A}(W)$ is a flat connection with $\onb(z)=Id_W$. Let
\[s:L\rightarrow L(f) \]
be defined as follows: $s(x):=(0,x)$. It follows $s$ is an $\B$-linear map satisfying
\[ s([x,y])=[s(x),s(y)]-f(x,y)z \]
for all $x,y\in L$. Define $\nabla:=\onb\circ s$. It follows we get a $\B$-linear map
\[ \nabla:L\rightarrow \End_{A}(W) \]
with the property that for all $a\in B, x\in L$ and $w\in W$ the following formula holds:
\[ \nabla(x)(aw)=a\nabla(x)(w)+\alpha(x)(a)w \]
hence $\nabla$ is an $L$-connection. We calculate the curvature $R_{\nabla}$: Let $x,y\in L$.
It follows
\[ R_{\nabla}(x,y):= \]
\[ [\nabla(x),\nabla(y)]-\nabla([x,y])=\]
\[ [\onb(s(x)),\onb(s(y))]-\onb(s([x,y]))=\]
\[ [\onb(s(x)),\onb(s(y))]-\onb([s(x),s(y)]-f(x,y)z)=\]
\[ R_{\onb}(s(x),s(y))+f(x,y)\onb(z)=f(x,y)Id_W \]
since $R_{\onb}=0$ by assumption. It follows $\nabla$ has curvature type $f$ and the Lemma follows.
\end{proof}

\begin{proposition} \label{fmodules}Let $W$ be a left $\B$-module and let $f\in \Z^2(L,\B)$.
There is a one-to-one correspondence between the set of left $U(B,L,f)$-module
structures on $W$ and the set of $L$-connections of curvature type $f$ on $W$.
\end{proposition}
\begin{proof} Recall the following: $L(f):=\B z\oplus L$ with $\alpha_f(az+x):=\alpha(x)$. Let
$\pb, \pl$ and $\plf$ be the maps defined above.
Let $W$ be a left $U(\B,L,f)$-module. Define for any $x\in L$ and $w\in W$ the following
map: $\nabla(x)(w):=\pl(x)w$. One checks that $\nabla$ is an $L$-connection on $W$. Assume $x,y\in L$ and $w\in W$.
It follows 
\[ \pl(x)\pl(y)-\pl(y)\pl(x)=\pl([x,y])+p_{\A}(f(x,y)) \]
in $U(\B,L,f)$ hence
\[   [\nabla(x),\nabla(y)](w)=\nabla([x,y])(w)+f(x,y)w .\]
It follows 
\[ R_\nabla(x,y)w=f(x,y)w\]
hence $\nabla$ is an $L$-connection of curvature type $f$.
Conversely let $\nabla:L\rightarrow \End_{\A}(W)$ be an $L$-connection of curvature type $f$.
Define the following action
\[\phi : T^1(L(f)) \rightarrow \End_{\A}(W) \]
by
\[\phi(\otimes_i (b_iz+x_i)):=\prod_i (b_iId_W+\nabla(x_i)).\]
One checks the action $\phi$ gives a map
\[ U(\B,L,f) \rightarrow \End_{\A}(W) .\]
One checks this construction sets up the desired correspondence and the Proposition is proved.
\end{proof}

Let $f\in \Z^2(L,\B)$ and recall the extension $L(f):=Bz\oplus L$ corresponding to $f$. Let 
$p:L(f)\rightarrow U(\B,L(f))$ be the canonical map where $U(\B,L(f))$ is Rineharts universal enveloping algebra
of $L(f)$. Let $z':=p(z)$ in $U(\B,L(f))$ and let $I:=(z'-1)\subseteq U(\B,L(f))$ be the 2-sided ideal generated by 
$z'-1$.  Let $U(\B,L(f),z')=U(\B,L(f))/I$.

\begin{proposition} \label{lconnf} Let $W$ be a left $\B$-module.
There is a one-to-one correspondence between the set of left $U(\B,L(f),z')$-module structures on $W$  and the set of
$L$-connections of curvature type $f$ on $W$.
\end{proposition}
\begin{proof} A left $U(\B,L(f),z')$-module $W$ corresponds to a flat $L(f)$-connection
\[ \nb:L(f)\rightarrow \End_{\A}(W)\]
with $\nb(z)=Id_W$. By Lemma \ref{conncurvf} it follows $\nb$ corresponds to an $L$-connection $\onb$ of curvature type $f$
and the Proposition follows.
\end{proof}

Assume $\alpha:L\rightarrow \Der_{\A}(\B)$ is a Lie-Rinehart algebra which is projective as left $\B$-module.
Assume $f\in \Z^2(L,\B)$ is a 2-cocycle of $L$. Let $\mod(L,f)$ be the category of $L$-connections of curvature type
$f$.

\begin{corollary} \label{modules} Let $(U,U_i)$ be a filtered algebra of type $\alpha$, where $\alpha:L\rightarrow \Der_R(A)$
is a Lie-Rinehart algebra and $U_1/U_0\cong L$. Assume $L$ is a projective $U_0$-module. Let $\mod(U)$ be the category 
of left $U$-modules. It follows there is an element $f\in \Z^2(L,\B)$  and an equivalence of categories
\[ \mod(U)\cong \mod(L,f) .\]
\end{corollary}
\begin{proof} The Corollary follows from Theorem \ref{maindeform} and Proposition \ref{fmodules}
since $U\cong U(\B,L,f)$ for some $f\in \Z^2(L,\B)$.
\end{proof}

\begin{example} \text{Modules on sheaves of rings of differential operators.}

Let $\pi:\Spec(\B)\rightarrow \Spec(\A)$ be an arbitrary morphism of affine schemes and let $X:=\Spec(\B), S:=\Spec(\A)$.
Let $\{\Gamma, \Gamma_i\}$ be
a filtered sheaf of associative rings on $X$ relative to $S$ satisfying axioms \ref{si1} - \ref{si5}. 
Assume $\H^0(X, \Gamma_1/\Gamma_0)$ is a projective $\H^0(X,\Gamma_0)$-module. It follows from Theorem
\ref{maindeform} there is a Lie-Rinehart algebra
\[ \alpha:L\rightarrow \Der_{\A}(\B) \]
together with an element $f\in \Z^2(L,\B)$ and an isomorphism $\H^0(X, \Gamma)\cong U(\B,L,f)$ of filtered associative rings.
The category of sheaves of left $\Gamma$-modules is by Corollary \ref{modules} equivalent to the category $\mod(L,f)$ of
$L$-connections of curvature type $f$.
Hence Lie-Rinehart algebras and $L$-connections arise naturally when studying deformations of
sheaves of rings of differential operators and categories of sheaves of modules on sheaves of rings of differential operators. 
\end{example}

In the following we give an interpretation of the morphisms in $\Adef (\Sym_{\B}^*(L))$ in terms of the Lie-Rinehart complex
$\C^p(L,\B)$. Let $f,g\in \Z^2(L,\B)$ and let $h\in \C^1(L,\B)$ with $d^1h=f-g$.
Let
\[ p_f:L(f)\rightarrow U(\B,L,f)\]
and
\[ p_g:L(g)\rightarrow U(\B,L,g) \]
be the canonical maps of left $\B$-modules.
Let
\[ \phi:T^1(L(f))\rightarrow U(\B,L,g) \]
be the following map:
\[ \phi((a_1z+x_1)\otimes \cdots \otimes (a_kz+x_k)):=p_g((a_1+h(x_1)z))\cdots p_g((a_k+h(x-k))z+x_k).\]
Let $u:=az+x, v:=bz+y \in L(f)$. It follows
\[ [u,v]=(x(b)-y(a)+f(x,y))z+[x,y],\]
where we write $x(b):=\alpha(x)(b)$.
We get the following calculation:
\[ \phi(u\otimes v-v\otimes u -[u,v])=\]
\[p_g((a+h(x))z+x)p_g((b+h(y))z+y)-p_g((b+h(y))z+y)p_g((a+h(x))z+x) \]
\[ -p_g((x(b)-y(a)+f(x,y)+h([x,y]))z+[x,y]).\]
In the following we drop writing $p_g$ since all calculations take place in the algebra $U(\B,L,g)$.
We get
\[ (a+h(x))(b+h(y))+(a+h(x))y+x(b+h(y))+xy-(b+h(y))(a+h(x))\]
\[ -(b+h(y))x-y(a+h(x))-yx -x(a)+y(b)-f(x,y)-h([x,y])-[x,y]=\]
\[ ay+h(x)y+x(b)+bx+x(h(y))+h(y)x+xy-bx-h(y)x-y(a)\]
\[ -ay-y(h(x))-h(x)y-yx -x(a)+y(b)-f(x,y)-h([x,y])-[x,y]=\]
\[ xy-yx+x(h(y))-y(h(x))-h([x,y])-[x,y]-f(x,y)=\]
\[xy-yx-g(x,y)-[x,y]=p_g(x)p_g(y)-p_g(y)p_g(x)-p_g([x,y])=0.\]
It follows $\phi$ descends to a map
\[ \phi':U(L(f))\rightarrow U(\B,L,g).\]
Recall the following 2-sided ideal in $U(L(f))$: 
\[J_f:=\{p_f(bw)-p_f(bz)p_f(w):b\in \B, w\in L(f)\}.\]
It follows
\[ \phi'(p_f(bw)-p_f(bz)p_f(w))=\]
\[\phi'(b(az+x)-bz\otimes (az+x))=\]
\[(ba+h(bx))+bx-b(a+h(x)+x)=h(bx)-bh(x)=0.\]
It follows $\phi'$ descends to a well defined map
\[ \theta_h:U(\B,L,f)\rightarrow U(\B,L,g) \]
defined by
\[ \theta_h(\prod_i (a_iz+x_i)):=\prod_i ((a_i+h(x_i))z+x_i).\]
We may construct a similar map
\[ \theta_{-h}:U(B,L,g)\rightarrow U(\B,L,f)\]
and one verifies that
\[ \theta_h\circ \theta_{-h}=\theta_{-h}\circ \theta_h=id.\]
Hence $\theta_h$ is an isomorphism with inverse $\theta_{-h}$.

Consider the canonical map
\[ p_f^k:L(f)^{\otimes k}\rightarrow U_k(\B,L,f), \]
where we write
\[ p_f^k(u_1\otimes \cdots \otimes u_k)=u_1u_2\cdots u_k.\]
It follows for any permutation $\sigma$ of $k$ elements the following formula holds:
\[ u_1\cdots u_k=u_{\sigma(1)}\cdots u_{\sigma(k)}+w\]
where $w\in U_{k-1}(\B,L,f)$. We get a well defined map
\[ P^k_f:\Sym_{\B}^k(L(f))\rightarrow U_k(\B,L,f)/U_{k-1}(\B,L,f) \]
inducing a canonical map
\[ \gamma_f^k:\Sym_{\B}^k(L)\rightarrow U_k(\B,L,f)/U_{k-1}(\B,L,f).\]
By definition it follows that for any element
\[ x_1\cdots x_k \in U_k(\B,L,f)/U_{k-1}(\B,L,f)\]
with $x_i \in L$, it follows
\[ \theta_h(x_1\cdots x_k)=(h(x_1)+x_1)\cdots (h(x_k)+x_k)\in U_k(\B,L,g)/U_{k-1}(\B,L,g).\]
It follows
\[ \theta_h(\gamma^k_f(x_1\cdots x_k))=\gamma^k_g(x_1\cdots x_k)+w\]
where $w\in U_{k-1}(B,L,g)$ hence
\[ \theta_h(\gamma^k_f(x_1\cdots x_k))=\gamma^k_g(x_1\cdots x_k) \]
for any $x_1\cdots x_k\in \Sym_B^*(L)$.

\begin{proposition} The map $\theta_h: U(\B,L,f) \rightarrow U(\B,L,g)$ is a map in $\Adef (\Sym_{\B}^*(L))$ with inverse
$\theta_{-h}$.
\end{proposition}
\begin{proof} 
By the above calculations we get a commutative diagram
\[
\diagram   Gr(U(\B,L,f)) \rto^{Gr(\theta_h)} & Gr(U(\B,L,g)) \\
           \Sym_{\B}^*(L) \uto^{\gamma_f} \urto^{\gamma_g} &  \\
\enddiagram.
\]
The Proposition follows.
\end{proof}

Since all maps in $\Adef (\Sym_{\B}^*(L))$ are isomorphisms if follows that if $\overline{f}\neq\overline{g}$ in $\H^2(L,\B)$
there are no maps $\theta:U(\B,L,f)\rightarrow U(\B,L,g)$ in $\Adef (\Sym_{\B}^*(L))$. 

Assume $d^1h=g-f$ with $f,g\in \Z^2(L,\B)$ and $h\in \C^1(L,\B)$ and assume
\[ \theta:U(\B,L,f)\rightarrow U(\B,L,g) \]
is a map in $\Adef (\Sym_{\B}^*(L))$.
Let $u_k(\B,L,f):=U_k(\B,L,f)/U_{k-1}(\B,L,f)$ and consider the canonical isomorphisms
\[ \gamma_f^k: \Sym_{\B}^k(L)\rightarrow u_k(\B,L,f) \]
and
\[ \gamma^k_g:\Sym_{\B}^k(L)\rightarrow u_k(\B,L,g).\]
Since $\theta$ is a map in $\Adef (\Sym_{\B}^*(L))$ it follows
\[  Gr(\theta)_k\circ \gamma^k_f =\gamma^k_g\]
hence
\[ Gr(\theta)_k=\gamma^k_g\circ (\gamma^k_f)^{-1}.\]
It follows we get a commutative diagram
\[
\diagram   \B \dto^{=} \rto & U_1(\B,L,f) \dto^{\theta} \rto & L  \dto^{=}   \\
           \B \rto     & U_1(\B,L,g) \rto                    & L 
\enddiagram.
\]
Since $\theta: U_1(\B,L,f)\rightarrow U_1(\B,L,g)$ is a map of Lie-Rinehart algebras and $L$ is a projective $\B$-module
it follows there are isomorphisms $U_1(\B,L,f)\cong L(f)$ and $U_1(\B,L,g)\cong L(g)$ as Lie-Rinehart algebras. 
We moreover get
\[ \theta(p_f(az+x))=p_g((a+h(x))z+x) \]
for some $h\in \C^1(L,\B)$ where $d^1h=g-f$.  Let $x_1,..,x_k\in L$.
We get
\[ \theta(p_f(x_1)\cdots p_f(x_k))=\theta(p_f(x_1))\cdots \theta(p_f(x_k))=\]
\[ p_g(h(x_1)z+x_1)\cdots p_g(h(x_k)z+x_k) \]
hence $\theta=\theta_h$ for some $h\in \C^1(L,\B):=\Hom_{\B}(L,\B)$.  Let $\Adef :=\Adef (\Sym_{\B}^*(L))$.
And let $U_f:=U(\B,L,f)$.

\begin{theorem} \label{morphisms} Let $\alpha:L\rightarrow \Der_R(A)$ be a Lie-Rinehart algebra 
where $L$ is a projective $\B$-module. Let $f,g\in \Z^2(L,\B)$. 
If $\overline{f}=\overline{g}$ in $\H^2(L,\B)$ it follows
\[\Hom_{\Adef }(U_f,U_g)=\Z^1(L,\B).\] 
If $\overline{f}\neq \overline{g}$ in $\H^2(L,\B)$ it follows
$\Hom_{\Adef }(U_f,U_g)=(0)$.
\end{theorem}
\begin{proof} Assume $h,k\in\C^1(L,\B)$ with $d^1h=d^1k=g-f$ It follows $k-h\in \Z^1(L,\B)$ hence $k=h+\beta$
with $\beta\in \Z^1(L,\B)$. This shows that $\Hom_{\Adef}(U_f,U_g)=\Z^1(L,\B)$. The Theorem follows.
\end{proof}

\begin{corollary} Assume $f\in \Z^2(L,\B)$ where $L$ is a projective $\B$-module. 
It follows $\Aut_{\Adef }(U_f)=\Z^1(L,\B)$
\end{corollary}
\begin{proof} The Corollary is immediate from Theorem \ref{morphisms}.
\end{proof}

\begin{example} \text{Moduli spaces of sheaves of rings of differential operators.}

Let $X:=\Spec(\B)$, $S:=\Spec(\A)$ and $\pi:X\rightarrow S$ be a morphism of schemes. 
Let also $\mathbb{V}(\Z^1(L,\B)):=\Spec(\Sym_{\A}(\Z^1(L,\B)^*)$. It follows from Theorem \ref{morphisms} that the set 
of morphisms in $\Adef (\Sym_{\B}^*(L))$ are parametrized by the points of the scheme $\mathbb{V}(\Z^1(L,\B))$. Hence
the objects and morphisms of the deformation groupoid $\Adef $ are parametrized by the schemes
$\mathbb{V}(\H^2(L,\B))$ and $\mathbb{V}(\Z^1(L,\B))$. One may ask for a more functorial formulation and proof of Theorem 
\ref{maindeform} and \ref{morphisms} and the existence of a groupoid scheme structure
\[ s,t: \mathbb{V}(\Z^1(L,\B))\rightarrow \mathbb{V}(\H^2(L,\B)) \]
on the cohomology of the Lie-Rinehart algebra $(L,\alpha)$. If such a gropoid scheme structure exists
it follows the corresponding stack quotient
\[ [\mathbb{V}(\Z^1(L,\B))\rightarrow \mathbb{V}(\H^2(L,\B))] \]
may be viewed as the moduli space of sheaves of rings of differential operators on the affine scheme 
$X$ relative to the base $S$. One may ask similar questions for the higher order cohomology groups $\H^i(L,\B)$.
\end{example}

\section{Application II: Algebraic connections on projective modules}

We use the constructions in the previous sections to study algebraic connections of curvature type
$f$ on finitely generated projective $\B$-modules.  We prove that for any Lie-Rinehart algebra $L$ that is
projective as $\B$-module and any cohomology class $c$ in $\H^2(L, \B)$, there is a finitely generated projective $\B$-module
$E$ with $c_1(E)=c$ (see Theorem \ref{main}). We also construct families of mutually non-isomorphic $\B$-modules 
of arbitrary high rank parametrized by $\H^2(L,\B)$. 
We prove for any Lie-Rinehart algebra $L$ which is projective as left $\B$-module, 
the existence of a subring $Char(L)$ of $\H^{2*}(L,\B)$. The ring $Char(L)$is non-trivial if and only if $\H^2(L,\B)\neq 0$.
We prove that $Char(L)$ is a  subring of the image of the Chern -character 
\[ Ch_{\Q}:\K(L)_{\Q} \rightarrow H^*(L,\B).\]
The definition of $Char(L)$ does not involve the Grothendieck group $\K(L)_{\Q}$. The problem of calculating generators for 
$\K(L)_{\Q}$ is an unsolved problem in general.

Assume $\A\rightarrow \B$ is a map of commutative rings where $\A$ contains a field $k$ of characteristic zero.
Let $E$ be a finitely generated projective left $\B$-module of rank $r$. 
There are two equivalent ways of defining a connection on 
$E$. Let $L:=\Der_{\A}(\B)$ and let $\Omega^1:=\Omega^1_{\B/\A}$ be the 
module of K\"{a}hler differentials. Let $e_1,\ldots ,e_r \subseteq E$
and $x_1,\ldots ,x_r \subseteq E^*$ be a \emph{projective basis} for $E$ in the sense of \cite{maa0}.
This means the following equation holds in $E$ for all $e\in E$:
\[ \sum_i x_i(e)e_i=e.\]
Note: A projective basis is sometimes referred to as a \emph{dual basis} (see \cite{anderson}).
One uses the projective basis $V=(e_i,x_j)$ to define connections
\[ \nabla:L\rightarrow \End_\A(E) \]
by
\[ \nabla(z)(e):=\sum_i z(x_i(e))e_i \]
and
\[ \overline{\nabla}:E\rightarrow \Omega\otimes_{\B} E\]
by
\[ \overline{\nabla}(e):=\sum_i d(x_i(e))\otimes e_i.\]
The curvature $R_{\nabla}$ of $\nabla$ is defined as 
\[ R_{\nabla}(x,y):= [\nabla(x),\nabla(y)]-\nabla([x,y])\]
for $x,y\in L$. We get a map
\[ R_{\nabla}:L\wedge_{\B} L \rightarrow \End_{\B}(E).\]
The connection $\overline{\nabla}$ gives rise to the algebraic DeRham complex
\[ E\rightarrow^{\nabla} \Omega^1 \otimes_{\B} E \rightarrow^{\nabla^1} E\otimes_{\B} \Omega^2 \rightarrow \cdots \]
and the curvature $K_{\overline{\nabla}}$ is defined as
\[ K_{\overline{\nabla}}:=\nabla^1 \circ \nabla .\]
Given the projective basis $V$ one gets an idempotent $\phi$ for the projective module $E$. The projective 
basis $V$ gives by the results in \cite{maa0} a surjection $p:\B^r \rightarrow E$ of left $\B$ modules with 
left $\B$-linear splitting $s:E\rightarrow \B^r$. The endomorphism $\phi=s\circ p\in \End_{\B}(\B^r)$ is an idempotent 
for $E$. In \cite{maa0} 
the following formula is proved for all $x,y\in L$:
\[ R_{\nabla}(x,y)= [x(\phi), y(\phi)]    \]
where $x(\phi)$ is the matrix we get when we let $x$ act on the coefficients of $\phi$. The product 
$[-,-]$ is the Lie product in the ring $\End_{\B}(\B^r)$.
Given a connection $\nabla: L\rightarrow \End_{\A}(E)$ all connections $\nabla'$ on $E$ are given as follows:
\[ \nabla'=\nabla+\psi \]
with $\psi\in \Hom_{\B}(L, \End_{\B}(E))$. The set of connections on $E$ is a \emph{torsor} on the abelian 
group $\Hom_{\B}(L,\End_{\B}(E))$. 
Hence it is a difficult problem to decide if a given module $E$ has a flat algebraic connection.
One has to study the set of connections
\[ \{ \nabla+\psi: \psi\in \Hom_{\B}(L,\End_{\B}(E)) \}\]
which is a large set in general.

If one is interested in the curvature of a connection it is more natural to use the language of Lie-Rinehart algebras 
because of the existence of the universal enveloping algebra $U(\B,L,f)$ for $f\in \Z^2(L,\B)$. 
Using $U(\B,L,f)$  we can give explicit
examples of a connection $(E,\nabla)$ where the curvature $R_{\nabla}$ satisfy certain properties.
For connections $\overline{\nabla}$ where the curvature $K_{\overline{\nabla}}$ is the composite
of two maps in the algebraic DeRham complex, there is no natural definition of an algebra or coalgebra
with properties similar to $U(\B,L,f)$. 

Hence Lie-Rinehart algebras appear naturally in the deformation 
theory of almost commutative rings, the theory of Chern classes and in various cohomology theories
as indicated in the introduction.

\begin{definition}
Let $\alpha:L\rightarrow \Der_{\A}(\B)$  be a Lie-Rinehart algebra which is a finitely generated
and projective $B$-module. Let $f\in Z^2(L,\B)$ be a 2-cocycle
and let $U^k(\B,L,f)$ be the descending filtration of $U(\B,L,f)$. 
Let for any $k\geq 1$ and $i\geq 1$ $V^{k,i}(\B,L,f):=U^k(\B,L,f)/U^{k+i}(\B,L,f)$.
\end{definition}

It follows $U^k(\B,L,f)$ is a filtration of two sided ideals in $U(\B,L,f)$.
By definition it follows $V^{k,i}(\B,L,f)$ is a left and right $U(\B,L,f)$ module for all $k,i\geq 1$.
Assume $rk(L)=l$ as projective $\B$-module. It follows by the results in the previous section that
$U^k(\B,L,f)$ and $V^{k,i}(\B,L,f)$ are projective $\B$-modules for all $k,i\geq 1$. Let $r(k,i,f):=rk(V^{k,i}(\B,L,f))$.

\begin{lemma} For all $k,i\geq 1$ the following formula holds:
\[ r(k,i,f) =\binom{l+k+i-1}{l}-\binom{l+k-1}{l}.\]
\end{lemma}
\begin{proof} Let $\O$ be the structure sheaf of $X:=\Spec(\B)$ and let $\mathcal{L}$ be the $\O$-module corresponding to
$L$. Let $\mathcal{U}_{L,f}$ be the left $\O$-module corresponding to $U(\B,L,f)$. It follows from Corollary \ref{pbw}
that there is an open subset $U$ in $X$ and a set of generators $s_1,..,s_l$ for $\mathcal{L}$ as free $\O(U)$-module
with the following property:
\[ \mathcal{U}_{L,f}(U)=\O(U)\{s_1^{p_1}\cdots s_l^{p_l}: p_i\geq 0\} \]
as free left $\O(U)$-module.
It follows $V^{k,i}(\B,L,f)(U)$ is described as follows:
\[ V^{k,i}(\B,L,f)(U)=\O(U)\{s_1^{p_1}\cdots s_l^{p_l}:l\leq \sum_j p_j < k+i\}\]
as free left $\O(U)$-module.
The Lemma follows.
\end{proof}

Note: The rank of $V^{k,i}(\B,L,f)$ is independent of choice of cocycle $f\in \Z^2(L,\B)$.

Since $V^{k,i}(\B,L,f)$ is a left $U(\B,L,f)$-module we get for all $k,i\geq 1$ algebraic connections
\[ \nabla: L\rightarrow \End_{\A}(V^{k,i}(\B,L,f)) \]
of curvature type $f$. Recall from Proposition \ref{fmodules} 
that this means that for any $x,y\in L$ and $w\in V^{k,i}(\B,L,f)$ it follows
\[ R_\nabla(x,y)(w)=f(x,y)w.\]
Let 
\begin{align}
&\label{FF} F:=\frac{1}{r(k,i,f)}f\in Z^2(L,\B). 
\end{align}
We get by Proposition \ref{fmodules} a connection
\[ \tilde{\nabla}:L\rightarrow \End_{\A}(V^{k,i}(\B,L, F))\]
of curvature type $F$. Let $c:=\overline{f}\in \H^2(L,\B)$.

\begin{theorem} \label{main}Let $\A$ contain a field of characteristic zero and let $\alpha:L\rightarrow \Der_R(A)$ 
be a Lie-Rinehart algebra which is finitely generated and projective as left $\B$-module of rank $l$. 
Let $F$ be the 2-cocycle defined in (\ref{FF}).
The following holds: 
\[ c_1(V^{k,i}(\B,L,F))= c \in \H^2(L,\B).\]
\end{theorem}
\begin{proof} By the results in \cite{maa1} we may construct the first Chern class of $V^{k,i}(\B,L, F)$
in $\H^2(L,B)$ by taking the trace of the curvature $R_{\tilde{\nabla}}$. It follows
\[ tr(R_{\tilde{\nabla}})=tr(F Id)=\frac{1}{r(k,i,f)}f tr(Id)=f.\]
Hence
\[c_1(V^{k,i}(\B,L,F))=\overline{f}=c \in \H^2(L,\B).\]
The Theorem is proved.
\end{proof}

\begin{corollary} \label{chernclass} The assumptions are as in Theorem \ref{main}.
Any cohomology class in $\H^2(L,\B)$ is the first Chern class 
of a finitely generated projective $\B$-module.
\end{corollary}
\begin{proof} The Corollary follows from Theorem \ref{main} since $f$ is an arbitrary element in $\Z^2(L,\B)$.
\end{proof}

\begin{example} \text{Holomorphic Lie algebroids.}

In the paper \cite{tortella} the following result is proved. Let $X$ be a smooth projective variety over the complex numbers
and let $\alpha:\mathcal{L}\rightarrow T_X$ be a holomorphic Lie algebroid. A holomorphic Lie algebroid is a complex analytic 
version of the notion of a Lie-Rinehart algebra. Let $P\in \H^0(X, \Omega^2_{\mathcal{L}})_{closed}$ 
and let $p$ be a numerical polynomial. If $P$ is not cohomologous to zero then Simpson's moduli space $M_{\mathcal{L},P}(p)$ 
from \cite{simpson} is empty.
Hence there is no locally free finite rank  $\O_X$-module $\mathcal{E}$ with a holomorphic $\mathcal{L}$-connection of 
curvature type $P$. In the affine situation
this is not true as Corollary \ref{chernclass} shows. Given a non-zero cohomology class $c\in \H^2(L,B)$ it follows
by the results in this section there are many non-isomorphic finitely generated projective modules $E$ with $c_1(E)=c$.
Note that a finite rank holomorphic vectorbundle on $X$ is algebraic hence $\mathcal{L}$ may be viewed as 
a \emph{sheaf of Lie-Rinehart algebras} on $X$.
\end{example}

\begin{corollary} Let $\alpha:L\rightarrow \Der_{\A}(\B)$ be a Lie-Rinehart algebra where $L$ is finitely generated and
projective as $\B$-module and let $X:=\Spec(\B)$. Let $M(X,L,c)$ be the moduli space of $L$-connections 
$(W,\nabla)$ where $W$ is a finitely generated
projective $\B$-module, $c\in \H^2(L,\B)$ and $c_1(W)=c$. It follows $M(X,L,c)$ is non-empty.
\end{corollary}
\begin{proof} The proof follows from the discussion above and Theorem \ref{main}.
\end{proof}

\begin{corollary} \label{chernsurjective} The assumptions are as in Theorem \ref{main}. The first Chern class map
\[ c_1:\K_0(L) \rightarrow \H^2(L,\B) \]
is a surjective map of abelian groups.
If $\H^2(L,B)\neq 0$ it follows $\K_0(L) \neq 0$.
\end{corollary}
\begin{proof} The map $c_1$ is a map of abelian groups. By Theorem \ref{main} is is surjective. 
If $0\neq c\in \H^2(L,B)$ it follows there is an $0\neq x\in \K_0(L)$ with $c_1(x)=c$ since $c_1$ is a map of
abelian groups. Hence $\K_0(L)$ is non-trivial. The Corollary follows. 
\end{proof}

\begin{example}\label{singularcohomology} \text{Singular cohomology of an affine algebraic variety.}

Let $A$ be a finitely generated regular algebra over the complex numbers and let $X:=\Spec(A)$ be the associated
affine scheme. Let $X_{\Cx}$ be the underlying complex algebraic manifold of $X$.
Corollary \ref{chernsurjective} gives a surjective map of abelian groups
\begin{align}
&\label{cplx}c_1:\K_0(\Der_{\Cx}(A))\rightarrow \H^2_{sing}(X_{\Cx},\Cx)
\end{align}
where $\H^2_{sing}(X_{\Cx},\Cx)$ is singular cohomology of $X_{\Cx}$ with complex coefficients. 
Since $\K_0(\Der_{\Cx}(A))$ is an abelian group there is a decomposition
\[ \K_0(\Der_{\Cx}(A))\cong T\oplus (\oplus_{i\in I}\mathbb{Z}e_i) \]
where $\oplus_{i\in I}\mathbb{Z}e_i$ is the free abelian group on the set $I$ and $T$ is a torsion group. 
If $\H^2_{sing}(X_{\Cx},\Cx)\neq 0$ and $T=0$ it follows $I$ is uncountable. Hence $\K_0(\Der_{\Cx}(A))$ may be large.

We get the following Corollary:
\begin{corollary} For any topological
class $c\in \H^2_{sing}(X_{\Cx},\Cx)$, there is a finite rank complex algebraic vector bundle $E$ on $X_{\Cx}$ with 
$c_1(E)=c$. 
\end{corollary}
\begin{proof} The Corollary follows since $c_1$ is a surjective map of abelian groups.
\end{proof}

The map $c_1$ is defined as follows: 
\[ c_1(\sum_i n_i[E_i, \nabla_i])=\sum_i n_i \overline{tr(R_{\nabla_i})}=\sum_i n_ic_1(E_i).\]
Here $n_i$ are integers. 

For any algebraic cycle $\omega=\sum_i n_i[V_i]$ on $X$ there is an associated topological cohomology class
\[ \gamma(\omega)\in \H^{*}_{sing}(X_{\Cx},\Cx).\]
Since the first Chern class $c_1(E)$ of a finitely generated and projective $A$-module $E$ is defined using 
an algebraic connection on $E$ it follows $c_1(E)$ lies in the image of $\gamma$ in $\H^2_{sing}(X_{\Cx},\Cx)$. It follows
any cohomology class $c$ in $\H^2_{sing}(X_{\Cx},\Cx)$ is the class of an algebraic cycle $\omega=\sum_i n_i[V_i]$ on $X$ 
with integral coefficients $n_i$. Hence the affine algebraic situation differs much from the projective situation.

Assume $E$ is a finitely generated projective $A$-module and let $\mathbb{V}(E^*)$ be the vector bundle associated to
$E$. Let $\mathbb{V}(E^*)_{\Cx}$ be the underlying complex algebraic manifold of $\mathbb{V}(E^*)$ in the strong topology.
Assume $\mathbb{V}(E^*)_{\Cx}$ has a Hermitian metric. It follows we may using Milnor's notes (see \cite{milnor}) to define
Chern classes 
\[  c_i(\mathbb{V}(E^*)_{\Cx})\in \H^{2i}_{sing}(X_{\Cx},\mathbb{Z}) \]
where $\H^{2i}_{sing}(X_{\Cx},\mathbb{Z})$ is singular cohomology of $X_{\Cx}$ with integer coefficients. We may look 
at the canonical image $W$ of $\H^{2}_{sing}(X_{\Cx},\mathbb{Z})$ in $\H^2_{sing}(X_{\Cx},\Cx)$. It follows
the Chern class $c_1(\mathbb{V}(E^*)_{\Cx})$ is in $W$, which is a strict sub-space of $\H^2_{sing}(X_{\Cx},\Cx)$.
Hence the Chern class $c_1(\mathbb{V}(E^*)_{\Cx})$ from \cite{milnor} 
differs from the Chern class $c_1(E)$ defined in \cite{maa1} using an algebraic connection
\[ \nabla: \Der_{\Cx}(A)\rightarrow \End_{\Cx}(E).\]
Assume $\mathbb{V}(E^*)_{\Cx}$ has a Hermitian metric.
Let $i:\H^2_{sing}(X_{\Cx},\mathbb{Z})\rightarrow \H^2_{sing}(X_{\Cx},\Cx)$ be the canonical map and
define $c_1^M(E):=i(c_1(\mathbb{V}(E^*)_{\Cx}))\in \H^2_{sing}(X_{\Cx},\Cx))$

\begin{definition}Let $E$ be a finitely generated and projective $A$ module with the property that
$\mathbb{V}(E^*)_{\Cx}$ has a Hermitian metric. Define
\[ c_1^{mot}(E):=c_1(E)-c_1^M(E)\in \H_{sing}^2(X_{\Cx},\Cx).\]
Let $c_1^{mot}(E)$ be the \emph{motivic Chern class} of $E$.
\end{definition}

By the above arguments it follows $c_1^{mot}(E)$ is a non-trivial characteristic class.
One wants to study the relationship between the class $c_1^{dev}(E)$ and the topologial and algebraic properties
of the projective module $E$.
\end{example}

\begin{corollary} \label{isomod} Assumptions are as in Theorem \ref{main}. 
Fix $k,i\geq 1$ and let $f_1,f_2\in \Z^2(L,\B)$.  
Assume $\overline{f_1}\neq \overline{f_2}$ in  $\H^2(L,\B)$. It follows $V^{k,i}(\B,L,f_1)$ and 
$V^{k,i}(\B,L,f_2)$ are non-isomorphic as left $\B$-modules.
\end{corollary} 
\begin{proof} Assume $V^{k,i}(\B,L,f_1)\cong V^{k,i}(\B,L,f_2)$ as left $\B$-modules.
Since $A$ has characteristic zero, it follows
\[ c_1(V^{k,i}(\B,L,f_2))=d\tilde{f_2}=d\tilde{f_1}=c_1(V^{k,i}(\B,L,f_1)\]
in $\H^2(L,\B)$ where $d=rk(V^{k,i}(\B,L,f_j)$. This leads to a contradiction and the Corollary follows.
\end{proof}

\begin{example} \label{family} \text{Families of finitely generated projective modules.}

Let $f,g\in \Z^2(L,\B)$ and
assume $\overline{f}=\overline{g}\in \H^2(L,\B)$. It follows there is an isomorphism $U(\B,L,f)\cong U(\B,L,g)$ of filtered
algebras. It follows for all $k\geq 1$ there is an isomorphism
\[ U^k(\B,L,f)\cong U^k(\B,L,g) \]
of left and right $\B$-modules hence $V^{k,i}(\B,L,f)\cong V^{k,i}(\B,L,g)$ as left and right $\B$-modules for all
$k,i\geq 1$. We may define for any cohomology class $c\in \H^2(L,\B)$
\[ V^{k,i}(\B,L,c):=V^{k,i}(\B,L,f) \]
where $f\in \Z^2(L,\B)$ is a representative for the class $c$.
Hence when we consider the left and right $\B$-module $V^{k,i}(\B,L,c)$ for varying $c \in \H^2(L,\B)$ 
we get a family of finitely generated projective $\B$-modules of constant rank 
parametrized by $\H^2(L,\B)$. 
From Lemma \ref{isomod} it follows that different classes in $\H^2(\B,L)$ gives non-isomorphic modules.
\end{example}

Recall that  $\A$ contains a field $k$ of characteristic zero and consider the map
\[ exp: \H^2(L,\B) \rightarrow \oplus_{k\geq 0} \H^{2k}(L,\B) \]
defined by
\[ exp(x):=\sum_{k\geq 0}\frac{1}{k!}x^k.\]

\begin{lemma} The map $exp$ is a map of abelian groups.
\end{lemma}
\begin{proof} We view the element $exp(x)$ as an element in the multiplicative subgroup of $\H^{2*}(L,\B)$
with ``constant term'' equal to one. 
Let $x,y\in \H^2(L,\B)$ be two cohomology classes.
We get
\[ exp(x+y):=\sum_{k\geq 0}\frac{1}{k!}(x+y)^k=\]
\[ \sum_{k\geq 0}\frac{1}{k!}\sum_{i+j=k}\binom{k}{i}x^iy^j=\]
\[\sum_{k\geq 0}\sum_{i+j=k}\frac{1}{i!j!}x^i y^j=(\sum_{i\geq 0}\frac{1}{i!}x^i)(\sum_{j\geq 0}\frac{1}{j!}y^j)=\]
\[exp(x)exp(y).\]
\end{proof}

Recall from \cite{maa2} the existence of a Chern-character
\[ Ch: \K(L)\rightarrow \H^{2*}(L,\B).\]
Extend $Ch$ to get a map
\[ Ch_{\mathbf{Q}}:\K(L)\otimes_{\mathbf{Z}}\mathbf{Q}\rightarrow \H^{2*}(L,\B).\]
Let $\nabla:L\rightarrow \End_{\A}(W)$ be an $L$-connection of curvature type $f$ with $f\in \Z^2(L,\B)$ 
and $W$ a finitely generated projective $\B$-module. Consider
the curvature $R_{\nabla}\in \C^2(L, \End_{\B}(W))$. We may use the shuffle product to get an element
\[ R_{\nabla}^k \in \C^{2k}(L,\End_{\B}(W)).\]
By definition
\[ Ch([W,\nabla]):=\sum_{k\geq 0} \frac{tr(R_{\nabla}^k)}{k!}\in \H^{2*}(L,\B).\]
Use the shuffle product to get the element $f^k\in \H^{2k}(L,\B)$. We get
\[f^k(x_1,\cdots ,x_{2k}):=\sum_{\sigma}(-)^{sgn(\sigma)}f(x_{\sigma(1)},x_{\sigma(2)})\cdots
f(x_{\sigma(2k-1)},x_{\sigma(2k)})\]
where the sum is over all $(2,2,..,2)$-shuffles.

\begin{lemma} The following formula holds for all $(x_1,..,x_{2k})\in L^{\times k}$:
\[ R_{\nabla}^k(x_1,..,x_{2k})=f^k(x_1,..,x_k)Id_W.\]
\end{lemma}
\begin{proof} We get the following calculation:
\[ R_{\nabla}^k(x_1,..,x_{2k}):=\sum_{\sigma}(-1)^{sgn(\sigma)}R_{\nabla}(x_{\sigma(1)},x_{\sigma(2)}\cdots
R_{\nabla}(x_{\sigma(2k-1)},x_{\sigma(2k)})=\]
\[ \sum_{\sigma}(-1)^{sgn(\sigma)}f(x_{\sigma(1)},x_{\sigma(2)}Id_W \cdots f(x_{\sigma(2k-1)},x_{\sigma(2k)})Id_W=\]
\[ f^k(x_1,..,x_{2k})Id_W.\]
The Lemma follows.
\end{proof}

It follows $tr(R_{\nabla}^k)=rk(W)f^k $.

\begin{definition} 
Let $\alpha:L\rightarrow \Der_R(A)$ be a Lie-Rinehart algebra whic is a finitely generated and projective $\B$-module and
assume $\H^2(L,\B)\neq 0$. Let $Char(L)$ be the subring of $\H^{2*}(L,\B)$ generated by the set
\begin{align}
&\label{setL} S:= \{ \sum r_i exp(x_i): x_i\in \H^2(L,\B), r_i\in \mathbf{Q} \}
\end{align}
Let $Char(L)$ be the \emph{characteristic ring of $L$}.
\end{definition}

The definition of the ring $Char(L)$ does not depend on a choice of a set of generators of $\K(L)$ since it is defined
in terms of $\H^2(L,\B)$. It is well known the problem of calculating generators of $\K(L)$ is an unsolved problem.

\begin{lemma} Let $\nabla:L\rightarrow \End_{\A}(E)$ and $\nabla':L\rightarrow \End_{\A}(F)$ be two connections with
\[ R_{\nabla}(x,y)=f(x,y)Id_E\]
and
\[ R_{\nabla'}(x,y)=g(x,y)Id_F \]
for $f,g\in \Z^2(L,\B)$. It follows
\[R_{\nabla \otimes \nabla'}(x,y)=(f(x,y)+g(x,y))Id_{E\otimes F}.\]
\end{lemma}
\begin{proof} One checks the module $E\otimes_{\B} F$ has a connection
\[ \eta(x)(u\otimes v)=\nabla(x)(u)\otimes v + u\otimes \nabla'(x)(v)\]
hence $\eta=\nabla \otimes \nabla'$.
It follows
\[ R_{\eta}(x,y)(u\otimes v)=R_{\nabla}(x,y)(u)\otimes v+u\otimes R_{\nabla'}(x,y)(v)=(f(x,y)+g(x,y))u\otimes v\]
and the Lemma follows.
\end{proof}

\begin{definition} Let $\overline{\K}(L)_{\Q}$ be the following sub-set of $\K(L)_{\Q}$:
\[ \overline{\K}(L)_{\Q}:=\{ \sum_i r_i[E_i,\nabla_i]:r_i \in \Q,R_{\nabla_i}=f_i, f_i\in \Z^2(L,B).\} \]
\end{definition}

\begin{proposition}\label{subchern} Let  $\alpha:L\rightarrow \Der_R(A)$ be a Lie-Rinehart algebra
which is finitely generated and projective as $\B$-module.
It follows the set $\overline{\K}(L)_{\Q}$ is a subring of $\K(L)_{\Q}$.
The Chern character $Ch$ induce a surjective map of rings
\[ \overline{Ch}:\overline{\K}(L)_{\Q}\rightarrow Char(L).\]
\end{proposition}
\begin{proof} Assume $x=\sum_i r_i[E_i,\nabla_i], y=\sum_j k_j[F_j,\nabla'_j]$ are elements in $\overline{\K}(L)_{\Q}$
It follows
\[ xy=\sum_{i,j}r_ik_j[E_i\otimes F_j, \nabla_i \otimes \nabla'_j] \]
and since
\[ R_{\nabla_i \otimes \nabla'_j}(x,y)=(f_i(x,y)+g_j(x,y))Id \]
it follows $xy\in \overline{\K}(L)_{\Q}$ hence $\overline{\K}(L)_{\Q}$ is closed under multiplication. The first claim
of the Proposition is proved.

Assume
\[ \sum_i r_i[E_i,\nabla_i]\in \overline{\K}(L)_{\Q}\]
with $R_{\nabla_i}(x,y)=f_i(x,y)Id$ It follows
\[ Ch(\sum_i r_i[E_i,\nabla_i])=\sum_i r_irk(E_i)exp(x_i) \]
where $x_i=\overline{f_i}\in \H^2(L,\B)$. It follows $Ch$ maps $\overline{\K}(L)_{\Q}$ into $Char(L)$. By definition
$Char(L)$ is the smallest subring of $\H^{2*}(L,\B)$ containing the set $S$ from \ref{setL}.
It follows $Char(L)$ is sums of products $s_1\cdots s_k$ where $s_i \in S$ for $i=1,..,k$ where $k\geq 0$.
Let
\[ s_j=\sum_i r^j_i exp(x^j_i)\in S\]
where $x^j_i=\overline{f^j_i}\in \H^2(L,B)$ where $f^j_i\in \Z^2(L,B)$. There is a connection $(E_{ji},\nabla_{ji})$
with $E_{ji}$ a finitely generated projective $\B$-module and where
\[R_{\nabla_{ji}}(x,y)=f^j_i(x,y)Id .\]
It follows
\[ Ch([E_{ji},\nabla_{ji}])=rk(E_{ji})exp(x^j_i).\]
We get
\[ Ch(\frac{1}{rk(E_{ji})}[E_{ji},\nabla_{ji}])=exp(x^j_i) \]
and
\[ Ch(\sum_i \frac{r^j_i}{rk(E_{ji})}[E_{ji},\nabla_{ji}])=\sum_i r^j_i exp(x^j_i)=s_j.\]
Let
\[ z_j=\sum_i \frac{r^j_i}{rk(E_{ji})}[E_{ji},\nabla_{ji}] .\]
Since $Ch$ is a ring homomorphism it follows
\[ Ch(z_1\cdots z_k)=Ch(z_1)\cdots Ch(z_k)=s_1\cdots s_k  .\]
It now follows  $\overline{Ch}$ is a surjective map. The Proposition is proved.
\end{proof}

\begin{corollary} Assumptions are as in \ref{subchern}. The characteristic ring 
$Char(L)$ is a sub ring of $Im(Ch)\subseteq \H^{2*}(L,\B)$. The ring $Char(L)$ is non-trivial if and only if
$\H^2(L,\B)\neq 0$.
\end{corollary}
\begin{proof} The Corollary follows from Proposition \ref{subchern} since $Ch$ is a map of rings.
\end{proof}

\begin{example} \text{An intrinsic description of the image of the Chern character.}

If $L$ is finitely generated and projective as $\B$-module and $\H^2(L,\B)\neq 0$ we get a non trivial extension
\[ 0\rightarrow Ker(\overline{Ch})\rightarrow \overline{\K}(L)_{\Q} \rightarrow Char(L) \rightarrow 0 \]
of rings. Hence $\K(L)_{\Q}$ is in a natural way an $\overline{\K}(L)_{\Q}$-module and 
$Im(Ch)$ is in a natural way a $Char(L)$-module. It is an unsolved problem to calculate 
generators for $\K(L)_{\Q}$ as $\overline{\K}(L)_{\Q}$-module and for $Im(Ch)$ as $Char(L)$-module.
One seek to construct a set $G=\{ E_i, \nabla_i\}_{i\in I}$ of connections
\[ \nabla_i:L\rightarrow \End_{\A}(E_i)\]
where $E_i$ is a finitely generated and projective $\B$-module for every $i\in I$,
with the property that the set $G$ generates $\K(L)_{\Q}$ as $\overline{\K}(L)_{\Q}$-module.
It then follows there is an equality
\[ Im(Ch)=Char(L)\{ Ch(x): x\in G\}\]
as subrings of $\H^{2*}(L,\B)$. 

Fix the notation of Example \ref{singularcohomology}. 
In the smooth projective case $Y\subseteq \mathbf{P}^n_{\Cx}$ the famous Hodge conjecture gives a conjectural description
of the rational span of the algebraic cycles as a sub-ring of $\H^*_{sing}(Y,\Cx)$ in terms of the Hodge decomposition. 
In the affine situation there is no Hodge-decomposition available. 
The characteristic ring $Char(\Der_{\Cx}(A))$ is a subring of $\H^*_{sing}(X_{\Cx},\Cx)$
consisting of algebraic cycles with rational coefficients. In the case when $\H^2_{sing}(X_{\Cx},\Cx)\neq 0$ and
$Char(\Der_{\Cx}(A))$ is large it may be the set $G$ is finite and one wants to give a formula for $G$ in terms of
invariants of $X$ in $\K_0(\Der_{\Cx}(A))$. One seek to give an intrinsic description of the image of the Chern character
in $\H^*_{sing}(X_{\Cx},\Cx)$ without choosing generators of the group $\K_0(\Der_{\Cx}(A))$. To calculate
generators for $\K_0(\Der_{\Cx}(A))$ is an unsolved problem in general.
\end{example}

\begin{example} \text{The Gauss-Manin connection.}

Let $\Cx\rightarrow \A \rightarrow \B$ be a sequence of maps of rings with $\A,\B$ finitely generated and regular over $\Cx$ .
Let $X:=\Spec(\B), S=:\Spec(\A)$ and $\pi:X\rightarrow S$ the induced morphism. Assume $\pi$ is smooth of relative dimension
$n$ hence $\Der_{\A}(\B)$ is a locally free $\B$-module of rank $n$. It follows we get an inclusion
\[ Char(\Der_{\A}(\B))\subseteq Im(Ch) \subseteq \H^{2*}(\Der_{\A}(\B), \B) .\]
There is an algebraic connection $\nabla_{GM}$ called the Gauss-Manin connection (see \cite{maa15})
\[ \nabla_{GM}: L \rightarrow \End_{\A}(\H^{2*}(\Der_{\A}(\B),\B)) .\]
We get for any $x$ in $L$ an endomorphism
\[ \nabla(x): \H^{2*}(\Der_{\A}(\B),\B)\rightarrow \H^{2*}(\Der_{\A}(\B),\B)  \]
acting as a differential operator of order one.
If for some $x$ the corresponding operator
 $\nabla(x)$ fix the subring $Im(Ch)$ we may use the connection 
$\nabla_{GM}$ and $Char(\Der_{\A}(\B))$ in the study
of $Im(Ch)$. One has to calculate explicitly the cohomology group $\H^i(\Der_{\A}(\B),\B)$ and the Gauss-Manin connection
$\nabla_{GM}$ for a class of smooth families $\pi$. Explicit formulas for algebraic connections 
have been calculated by hand in \cite{maa0} for a class of cotangent bundles on ellipsoid surfaces using the notion
of a projective basis. If $\H^{2*}(\Der_{\A}(\B),\B)$ is a finitely generated and projective $\A$-module, the techniques
from \cite{maa0} give explicit formulas for such connections.
\end{example}

\section{Appendix A: Categories of $L$-connections and module categories}

In this section we study the category of connections $\nabla:L\rightarrow \End_{\A}(W)$ where  $W$ is
finitely generated and projective as left $\B$-module, denoted $\conn^{fp}(L)$. We give an explicit realization of
 the category $\conn^{fp}(L)$ as a category of modules on an associative ring $\U(L)$ (see Corollary \ref{equiv}).
We use the associative ring $\U(L)$ to define $\Ext$ and $\Tor$-groups
\[ \Ext^i_{\U(L)}(V,W) , \Tor^i_{\U(L)}(V,W) \]
for any pair of $L$-connections $V,W$ and any integer $i\geq 0$ (see \ref{cohomology}).
A definition of $\Ext$ and $\Tor$-groups of $L$-connections was previously only given in the case of flat $L$-connection.

Assume in the following that $\A\rightarrow \B$ is a unital map of commutative unital rings and $\alpha:L\rightarrow \Der_R(A)$
an arbitrary Lie-Rinehart algebra.

\begin{definition} Let $W$ be an abelian group with a left and right $\B$-module structure.
If the following equation holds
\[ a(wb)=(aw)b \]
for all $a,b \in \B$ and $w\in W$ we say that $W $ is a $(\B,\B)$-module.
Assume $V,W$ are $(\B,\B)$-modules. A map of abelian groups
\[ \phi: V\rightarrow W\]
that is left and right $\B$-linear separately is a \emph{map of $(\B,\B)$-modules}.
\end{definition}

Assume
$\alpha:L\rightarrow \Der_{\A}(\B)$ is a Lie-Rinehart algebra. Let $\B z$ be the free rank one $\B$-module on the element
$z$ and consider the direct sum $\tilde{L}:=\B z\oplus L$. Define the following left and right actions on $\tL$:
Let
\begin{align}
&\label{m1} a(bz+x):=(ab)z+ax
\end{align}
and
\begin{align}
&\label{m2} (bz+x)a:=(ba+x(a))z+ax
\end{align}
for any element $a\in \B$. We write $\alpha(x)(a):=x(a) $ for simplicity.

\begin{lemma} \label{leftright} The actions (\ref{m1}) and (\ref{m2}) define an $(\B,\B)$-module structure on $\tL$.
Assume $L$ is finitely generated and projective as left $\B$-module. It follows $\tL$ is finitely generated
and projective as left and right $\B$-module separately. Moreover $\tL$ is a left $\B\otimes_{\A}\B$-module.
\end{lemma}
\begin{proof} The action (\ref{m1}) is clearly a left $\B$-module structure on $\tL$. One checks that (\ref{m2}) is a
right $\B$-module structure on $\tL$. One finally checks that for any $a,b \in \B$ and $w\in \tL$ the following holds:
\[ a(wb)=(aw)b .\]
Obviously $\tL$ is finitely generated and projective as left $\B$-module. We prove this statement also holds
when we view $\tL$ as right $\B$-module. Assume $p:\B^m\rightarrow L$ is a surjective map of left $\B$-modules with a
left $\B$-linear section $s:L\rightarrow \B^m$. Define the following maps:
\[ \tp:\B\oplus \B^m\rightarrow \tL \]
by
\[ \tp(a,u):=(a,p(u)) \]
and
\[ \ts:\tL\rightarrow \B\oplus \B^m \]
by
\[ \ts(b,x):=(b,s(x)).\]
It follows $\ts$ is a section of $\tp$. The module $\B\oplus \B^m$ is in a canonical way a left $\B$-module.
Define the following right action of $\B$ on $\B\oplus \B^m$:
\[ (b,u)a:=(ba+\alpha(p(u))(a), au)\]
for $a\in \B$. One checks that $\B\oplus \B^m$ becomes a left and right $\B$-module with this right action.
One moreover checks that $\tp$ and $\ts$ are left and right $\B$-linear. It follows that the map
\[ \tp:\B\oplus \B^m \rightarrow \tL \]
is right $\B$-linear with right $\B$-linear splitting given by $\ts$. It follows $\tL$ is finitely generated and
projective as right $\B$-module.
Define for any element  $a\otimes b\in \B\otimes_{\A} \B$ the following action:
\[ (a\otimes b).(bz+x):=a((bz+x)b).\]
One checks this gives a left action of $\B\otimes_{\A} \B$ on $\tL$.
\end{proof}

\begin{definition} Let $W$ be a left $\B$-module and let $\phi\in \End_{\B}(W)$. An $\B$-linear map
\[ \nabla:L\rightarrow \End_{\A}(W) \]
satisfying 
\[ \nabla(x)(aw)=a\nabla(x)(w)+\alpha(x)(a)\phi(w) \]
for $a\in \B$, $x\in L$ and $w\in W$ is an \emph{$(L,\phi )$-connection} on $W$. Let $\conn^{\End}(L)$ denote 
the cateory of $(L,\phi)$-connections for varying $\phi \in \End_{\B}(W)$ and left $\B$-modules $W$.
Let $\conn(L)$ denote the category of $(L, Id)$-connections. Let $\conn^{fp}(L)$ denote the category
of $L$-connections $(W,\nabla)$ where $W$ is a finitely generated and projective $\B$-module.
\end{definition}

Note: An $(L,Id_W)$-connection $\nabla:L\rightarrow \End_{\A}(W)$ is an ordinary $L$-connection.

If $V,W$ are $(\B,\B)$-modules, let $\Hom_{(\B,\B)}(V,W)$ denote the abelian group of $(\B,\B)$-linear
maps from $V$ to $W$. It follows $\Hom_{(\B,\B)}(V,W)$ is a $(\B,\B)$-module.

\begin{proposition} There is an equality of sets between the set $\Hom_{(\B,\B)}(\tL, \End_{\A}(W))$ and the set
of $(L,\phi)$-connections on $W$ for varying $\phi \in \End_{\B}(W)$.
\end{proposition}
\begin{proof}Assume 
\[ \rho: \tL\rightarrow \End_{\A}(W) \]
is an $(\B,\B)$-linear map. Let $u:=bz+x\in \tL$ and $a\in \B$ be elements.
It follows
\[ ua=(bz+x)a=(ba+x(a))z+ax=a(bz+x)+x(a)z =au+x(a)z \in \tL.\]
Consider the endomorphism $\rho(z)\in \End_{\A}(W)$. We get for $b\in \B$ and $w\in W$ the following calculation
\[ \rho(z)(bw)=(\rho(z)b)(w)=\rho(zb)(w)=\rho(bz)(w)=b\rho(z)(w),\]
hence $\phi= \rho(z) \in\End_{\B}(W)$.
We get for $a\in \B$
\[\rho(u)a=\rho(ua)=\rho(a(bz+x)+x(a)z)=a\rho(u)+x(a)\rho(z)=a\rho(u)+x(a)\phi.\]
It follows for $x\in L$, $a\in \B$ and $w\in W$ that
\[ \rho(x)(aw)=a\rho(x)(w)+x(a)\phi(w).\]
Hence the induced map
\[\rho:L\rightarrow \End_{\A}(W) \]
is a $\phi$-connection on $W$.

Conversely assume $\nabla:L\rightarrow \End_{\A}(W)$ is an $(L,\phi)$-connection for $\phi\in \End_{\B}(W)$.
Define the following map
\[\rho:\tL\rightarrow \End_{\A}(W)\]
by
\[ \rho(bz+x)(w):=b\phi(w)+\nabla(x)(w).\]
One checks that $\rho$ is a left $\B$-linear map.
We prove it is right $\B$-linear.
We get
\[ \rho((bz+x)a)(w)=\rho((ba+x(a))z+ax)(w)=(ba+x(a))\phi(w)+\nabla(ax)(w)=\]
\[ba\phi(w)+a\nabla(x)(w)+x(a)\phi(w)=a(b\phi(w)+\nabla(x)(w))+x(a)\phi(w)=\]
\[ a\rho(bz+x)(w)+x(a)\phi(w).\]
We get
\[\rho(bz+x)(aw)=b\phi(aw)+\nabla(x)(aw)=ab\phi(w)+a\nabla(x)(w)+x(a)\phi(w)=\]
\[a(b\phi(w)+\nabla(x)(w))+x(a)\phi(w)=a\rho(bz+x)(w)+x(a)\phi(w).\]
It follows
\[ \rho((bz+x)a)(w)=\rho(bz+x)(aw)\]
hence
\[\rho((bz+x)a)=\rho(bz+x)a.\]
It follows 
\[\rho:\tL\rightarrow \End_{\A}(W)\]
is an $(\B,\B)$-linear map. One checks this construction defines an equality of sets and the Proposition is proved.
\end{proof}

Let $\T_{\B\otimes_{\A} \B}(\tL)$ be the tensor algebra of $\tL$ as left $\B\otimes_{\A} \B$-module.

\begin{lemma} There is an isomorphism of abelian groups
\[\Hom_{(\B,\B)}(\tL, \End_{\A}(W)) \cong \Hom_{\B\otimes_{\A} \B-alg}(\T_{\B\otimes_{\A} \B}(\tL), \End_{\A}(W)) .\]
\end{lemma}
\begin{proof}  Assume 
\[ \phi:\tL\rightarrow \End_{\A}(W) \]
is a $(\B,\B)$-linear map. There are obvious $\B\otimes_{\A} \B$-actions on $\tL$ and $\End_{\A}(W)$ and we get
\[ \phi((a\otimes b).u)=\phi(a(ub))=a\phi(ub)=a(\phi(u)b)=(a\otimes b).\phi(u) \]
hence $\phi$ is $\B\otimes_{\A} \B$-linear. 
Because of the universal property of the tensor algebra of a module, it follows there is an isomorphism 
of abelian groups
\[ \Hom_{\B\otimes_{\A} \B}(\tL, \End_{\A}(W)) \cong \Hom_{\B\otimes_{\A} \B-alg}(\T_{\B\otimes_{\A} \B}(\tL), \End_{\A}(W)) .\]
The Lemma follows.
\end{proof}

Let $\mod(\T_{\B\otimes_{\A} \B}(\tL))$ denote the category of left $\T_{\B\otimes_{\A} \B}(\tL)$-modules.
Let $z' $ be the image of $z$ in the tensor algebra $\T_{\B\otimes_{\A} \B}(\tL)$ and let $I$ be the two sided ideal
generated by the element $z'-1$. Let $\U(L)$ be the quotient algebra $\T_{\B\otimes_{\A} \B}(\tL)/I$. Since there is an
equality between the category of $\B\otimes_{\A} \B$-algebra morphisms
\[ \rho: \T_{\B\otimes_{\A} \B}(\tL) \rightarrow \End_{\A}(W) \]
and the category $\mod(\T_{\B\otimes_{\A} \B}(\tL))$ it follows the category of $(L,\phi)$-connections
equals the category $\mod(\T_{\B\otimes_{\A} \B}(\tL))$.

\begin{lemma} \label{center} The center of $\U(L)$  contains the ring $A$.
\end{lemma}
\begin{proof} Let $w:=bz+x\in \tL$ and let $a\in \B$. It follows 
\[ wa=(ba+x(a))z+ax=x(a)z+a(bz+x)=aw+x(a)z.\]
Let $w_1,w_2\in \tL$ and consider the element $w_1\otimes w_2\in \tL\otimes_{\B\otimes_{\A} \B}\tL$.  It follows
\[ (w_1\otimes w_2)a=w_1\otimes (w_2a)=w_1\otimes (aw_2+x(a)z)=\]
\[ a(w_1\otimes w_2)+ w_1\otimes x(a)z .\]
It follows 
\[ (w_1\otimes w_2)a=a(w_1\otimes w_2)\]
if $a\in \A$ since $x(a)=0$. A similar calculation shows that for any element
\[ w_1\otimes \cdots \otimes w_k \in T_{\B\otimes_{\A} \B}(\tL)\]
and any $a\in \A$ it follows 
\[ (w_1\otimes \cdots \otimes w_k)a=a(w_1\otimes \cdots \otimes w_k).\]
Hence $\A$ lies in the center of $T_{\B\otimes_{\A} \B}(\tL)$. It follows $\A$ lies in the center of $\U(L)$ and the Lemma
follows.
\end{proof}

\begin{definition} Let $\U(L)$ be the \emph{universal algebra} for the category $\conn(L)$.
\end{definition}

The universal algebra $\U(L)$ differs from the universal enveloping algebra $U(\B,L,f)$ for $f\in \Z^2(L,\B)$
since a left $\U(L)$-module can be an arbitrary $L$-connection
\[ \nabla:L\rightarrow \End_{\A}(W).\]
It follows that a left $\U(L)$-module $W$ which is finitely 
generated and projective as $\B$-module
in a canonical way has an $L$-connection $\nabla$. Note that a left $\U(L)$-module $W$ in a is a $\B\otimes_{\A} \B$-module.
One checks that this is defined as follows: $(a\otimes b).w=(ab)w$. Hence $W$ is trivially an $(\B,\B)$-module.

Let $L(f)$ be the $\A$-Lie algebra with the Lie product defined as follows:
\[ [az+x, bz+y]:=(x(b)-y(a)+f(x,y))z+[x,y].\]
There is a canonical map of $\A$-Lie algebras
\[ p:L(f)\rightarrow U(\B,L,f).\]
Recall the structure of $(\B,\B)$-module on $L(f)$ defined as follows. Let $a\in \B$ and $bz+x\in L(f)$.
Define
\[ a(bz+x):=(ab)z+ax \]
and
\[ (bz+x)a:=(ba+x(a))z+ax.\]
It follows from Lemma \ref{leftright} $L(f)$ is an $(\B,\B)$-module.

\begin{lemma} \label{2ideal} The map $p$ is a map of $(\B,\B)$-modules inducing a surjective map 
\[ \tilde{p}:\U(L)\rightarrow U(\B,L,f) \]
of associative rings.
\end{lemma}
\begin{proof} Let $w:=bz+x\in L(f)$ and let $a\in B$. The map $p$ is clearly a linear map.
We get
\[p(aw):=p((ab)z+ax)=(ab)1+ax=a(b1+x)=ap(w) \]
hence $p$ is left $\B$-linear. We get
\[ p(wa):=p((ba+x(a))z+ax)=(ba+x(a))1+ax=(b1+x)a=p(w)a\]
hence the map $p$ is right $\B$-linear. 
Since the tensor algebra $\T_{\B\otimes_{\A} \B}(L(f))$ is universal with respect to maps
of $(\B,\B)$-modules it follows we get a canonical map
\[ p': \T_{\B\otimes_{\A} \B}(L(f))\rightarrow U(\B,L,f) \]
of associative rings. Since $p'(z-1)=1-1=0$ it follows $p'$ induce a map
\[ \tilde{p}:\U(L)\rightarrow U(\B,L,f) \]
of associative rings. The map $\tilde{p}$ is surjective by Proposition \ref{generated}. The Lemma follows.
\end{proof}

\begin{example}\label{twosided} \text{2-sided ideals in $\U(L)$ and curvature of $L$-connections.}

We get from Lemma \ref{2ideal} a canonical short exact sequence of associative rings
\begin{align}
&\label{uquot} 0\rightarrow I_{L,f}\rightarrow \U(L) \rightarrow U(\B,L,f)\rightarrow 0
\end{align}
for every $f\in \Z^2(L,B)$.
Hence the algebra $\U(L)$ is universal for the algebras $U(\B,L,f)$: For any element $f\in \Z^2(L,B)$ there is 
canonical a 2-sided
ideal $I_{L,f}$ in $\U(L)$ with $\U(L)/I_{L,f}\cong U(\B,L,f)$. It follows a left $\U(L)$-module $W$ which is annihilated
by the ideal $I_{L,f}$ corresponds to an $L$-connection $\nabla$ of curvature type $f$. We get a correspondence between
2-sided ideals $I$ in $\U(L)$ and $L$-connections $(W, \nabla)$, where the fact that the ideal $I$ annihilates $W$ imposes
restrictions on the curvature $R_\nabla$: If $I_{L,f}W=0$ it follows $R_\nabla=fId_w$.

Recall the definition of the decending sequence $U^k(\B,L,f)$ in $U(\B,L,f)$. We get from the exact sequence (\ref{uquot})
an exact sequence of left $\U(L)$-modules
\[0\rightarrow I_{L,f}\rightarrow U^k(L,f) \rightarrow U^k(\B,L,f)\rightarrow 0.\]
We get a filtration of left $\U(L)$-modules in $\U(L)$
\[ I_{L,f}\subseteq \cdots \subseteq U^2(L,f)\subseteq U^1(L,f)=\U(L)\]
with the property that there is an isomorphsm of $\U(L)$-modules
\[ U^k(L,f)/U^{k+i}(L,f)\cong V^{k,i}(\B,L,f).\]
Hence the $L$-connection
\[\nabla:L\rightarrow \End_{\A}(V^{k,i}(\B,L,f)) \]
may be constructed using filtrations in $\U(L)$. One may try to generalize this construction and construct
descending filtrations $\{U^i\}$ in $\U(L)$ whose successive quotients $U^i/U^{i+j}$ give rise to $L$-connections 
$(W_{i,j},\nabla_{i,j})$ where $W_{i,j}$ is finitely generated and projective as left $\B$-module. 
Properties of the filtration $(U^i)$ will impose restrictions on the curvature $R_{\nabla_{i,j}}$ of the 
connection $\nabla_{i,j}$. This will be applied to the study of the image of the Chern character.

We moreover get for any $f$ in $\Z^2(L,\B)$ exact sequences
\[ 0\rightarrow I_{L,f}\rightarrow U_k(L,f)\rightarrow U_k(\B,L,f)\rightarrow 0\]
for any integer $k\geq 0$. Let $U_{-1}(L,f)=I_{L,f}$. If $L$ is projective as left $\B$-module it follows 
there is for every $f\in \Z^2(L,\B)$ an isomorphism
\[ \Sym_{\B}^*(L)\cong \oplus_{k\geq 0}U_k(L,f)/U_{k-1}(L,f)\]
of graded $\B$-algebras.
Hence the universal algebra $\U(L)$ contains ascending and descending filtrations with properties similar to
$\{U_k(\B,L,f)\}_{k\geq 1}$ and $\{U^k(\B,L,f)\}_{k\geq 0}$ for all $f\in \Z^2(L,\B)$.
\end{example}

Let $\underline{Lie}_{\B/\A}$ be the category of Lie-Rinehart algebras $\alpha:L\rightarrow \Der_{\A}(\B)$ and maps
and let $\underline{Rings}$ be the category of associative rings with unit and unital maps of rings.

\begin{lemma} There is a covariant functor
\[\U:\underline{Lie}_{\B/\A}\rightarrow \underline{Rings}.\]
\end{lemma}
\begin{proof} Let $\U(L)$ be the universal algebra of $L$. One checks this defines a functor and the Lemma follows.
\end{proof}

\begin{definition}
Let $\mod(\U(L))$ denote the category of left $\U(L)$-modules. $\mod^{fp}(\U(L))$ be the category of left $\U(L)$-modules
which are finitely generated and projective as left $\B$-module. 
\end{definition}

\begin{example} \text{Projective objects in $\conn(L)$.}

Note: In the category $\conn(L)$ there is a definition of the notion of an injective and projective object. Recall the
definition of a projective object in $\conn(L)$.
An
$L$-connection $(P,\nabla)$ is \emph{projective} if and only if for any map of $L$-connections
\[ f:(P,\nabla)\rightarrow (W,\nabla_W) \]
and any surjective map of $L$-connections
\[u:(V,\nabla_V)\rightarrow (W,\nabla_W) \]
there is a map of connections
\[\tilde{f}:(P,\nabla)\rightarrow (V,\nabla_V) \]
with $u\circ \tilde{f}=f$. If $P$ is a projective $\B$ module it follows there is a lifting $\tilde{f}$ of $f$
as maps of $\B$-modules. There are in general many such liftings but these liftings are in general not maps of connections.
Hence it is not clear
how to construct an $L$-connection $(P, \nabla)$ that is a projective object in $\conn(L)$.
The category $\conn(L)$ is a small abelian category hence a well known result from category theory
says $\conn(L)$ may be realized as a sub category of the category $\mod(R)$ of left modules on an associative
ring $R$. This construction does not preserve injective and projective objects in general.
\end{example}

\begin{theorem} \label{equiv} There is an exact equivalence of categories
\begin{align}
&\label{morita}  \conn(L)\cong \mod(\U(L)) .
\end{align}
The equivalence preserve injective and projective objects.
\end{theorem}
\begin{proof} By the above argument it follows there is a one to one correspondence between $L$-connections
\[ \nabla:L\rightarrow \End_{\A}(V)\]
and left $\U(L)$-modules. Assume $(V,\nabla)$and $(W,\rho)$ are $L$-connection. There is a connection
\[ \eta:L\rightarrow \End_{\A}(\Hom_{\B}(V,W)) \]
defined as follows:
\[ \eta(x)(\phi):=\rho(x)\circ \phi -\phi\circ \nabla(x).\]
We define
\[ \Hom_{\conn(L)}(V,W):=\Hom_{\B}(V,W)^{\eta}, \]
where 
\[ \Hom_{\B}(V,W)^{\eta}:=\{ \phi:\eta(x)(\phi)=0\text{ for all }x\in L\}.\]
Let $\phi\in \Hom_{\U(L)}(V,W)$. It follows $\phi:V\rightarrow W$ is a $\B$-linear map.
Since $\phi(zv)=z\phi(v)$ for all $z\in \U(L)$ it follows $\phi(xv)=x\phi(v)$ for all $x\in L$. It follows
$x\in \Hom_{\B}(V,W)^{\eta}$. Assume $\phi\in \Hom_{\B}(V,W)^{\eta}$. It follows $\phi(xv)=x\phi(v)$ for all $x\in L$.
Let $u:=u_1\cdots u_k=(a_1z+x_1)\cdots a_kz+x_k)\in \U(L)$. We get
\[ \phi(u_1v)=\phi((a_1z+x_1)v)=\phi(a_1v+x_1v)=a_1\phi(v)+x_1\phi(v)=(a_1z+x_1)\phi(v)=u_1\phi(v).\]
We get
\[ \phi(uv)=\phi(u_1\cdots u_kv)=u_1\cdots u_k\phi(v)=u\phi(v) \]
by induction. It follows $\phi\in \Hom_{\U(L)}(V,W)$. It follows we get an exact equivalence
of categories. It is clear this equivalence preserve injective and projective objects. The Theorem is proved.
\end{proof}

There is work in progress giving a global version of Theorem \ref{equiv} valid for any sheaf of Lie-Rinehart algebras
on an arbitrary scheme $X$.

\begin{example} \text{Projective objects in $\conn(L)$ II.}

Theorem \ref{equiv} gives an elementary, explicit and functorial
construction of an associative ring $\U(L)$ and an equivalence between $\mod(\U(L))$ and $\conn(L)$. This equivalence
preserves projective and injective objects. Hence an $L$-connection $(P, \nabla)$ is a projective (resp. injective) object
in $\conn(L)$ if and only if the corresponding left $\U(L)$-module $P$ is a projective (resp. injective) $\U(L)$-module.
Any left $\U(L)$-module $P$ has a surjection
\[ p:\oplus_{i\in I } \U(L)e_i \rightarrow P \rightarrow 0\]
of left $\U(L)$-modules, where $I$ is some index set. The following result is well known.

\begin{lemma} Assume $(P,\nabla)$ is an $L$-connection. It follows $P$ is projective as left $\U(L)$-module 
if and only if $P$ is a direct summand of a free
$\U(L)$-module.
\end{lemma}
\begin{proof} The proof is an exercise.
\end{proof}

Hence projective objects in $\conn(L)$ are $L$-connections $(P,\nabla)$ where the underlying
$\U(L)$-module of $P$ is a direct summand of a free $\U(L)$-module.
\end{example}

\begin{example} The generalized Atiyah sequence and generalized Atiyah class.\end{example}

In this example we define the generalized Atiyah sequence for a Lie-Rinehart algebra $\alpha:L\rightarrow \Der_{\A}(\B)$
in the relative situation.
The generalized Atiyah sequence gives for any left $\B$-module $E$ rise to the generalized Atiyah class
\[ a_{L}(E)\in \Ext^1_{\B\otimes_{\A}\B}(L\otimes_{\A} E, E).\]
The generalized Atiyah class is zero if and only if $E$ has an $L$-connection
\[ \nabla:L\rightarrow \End_{\A}(E).\]

Let $\alpha:L\rightarrow \Der_{\A}(\B)$ be a Lie-Rinehart algebra and let $E$ be an arbitrary left $\B$-module.
Consider the $\A$-module $L\otimes_{\A} E$. The left $\B$-modules $L,E$ has a canonical right $\B$-module structure
defined by $xa:=ax$ for $x\in L$ and $a\in \B$. Moreover $eb:=be$ for $e\in E$. It follows $L\otimes_{\A} E$ has a canonical
left and right $\B$-module structure with the property that for any $w\in L\otimes_{\A} E$ it follows
$a(wb)=(aw)b$ for any elements $a,b\in \B$. It follows $L\otimes_{\A} E$ is an $(\B,\B)$-module and a $\B\otimes_{\A} \B$-module.

\begin{definition} An $\A$-linear map
\[ D:L\otimes_{\A} E\rightarrow E\]
is $(\B,\Der)$-linear if the following holds for all $a\in \B$ and $x\otimes e\in L\otimes_{\A} E$:
\[ D(a(x\otimes e))=aD(x\otimes e)\]
and
\[D((x\otimes e)a)=D(x\otimes e)a+\alpha(x)(a)e.\]
We say $D$ \emph{is a left $\B$-linear map and a right derivation}.
\end{definition}

\begin{lemma}\label{onetoone} There is a one to one correspondence between the set of $(\B,\Der)$-linear maps
\[ D:L\otimes_{\A} E\rightarrow E\]
and the set of $L$-connections
\[ \nabla:L\rightarrow \End_{\A}(E).\]
\end{lemma}
\begin{proof} Assume $\nabla:L\rightarrow \End_{\A}(E)$ is an $L$-connection.
Define the following map
\[ D:L\otimes_{\A} E \rightarrow E\]
by
\[ D(x\otimes e):=\nabla(x)(e).\] One checks that $D$ is a well defined $(\B,\Der)$-linear map.
Assume conversely that $D:L\otimes_{\A} E\rightarrow E$ is a $(\B,\Der)$-linear map.
Define
\[\nabla:L\rightarrow \End_{\A}(E) \]
by $\nabla(x)(e):=D(x\otimes e)$. It follows $\nabla$ is a well defined connection. The Lemma follows.
\end{proof}

Let $J^1_L(E):=E\oplus L\otimes_{\A} E$ be the \emph{first order jet module of type $L$} of $E$. It follows
$J^1_L(E)$  is a left $\B$-module. Define the following right action of $\B$ on $J^1_L(E)$:
\[ (e,x\otimes f)a:=(ea+x(a)f, x\otimes (fa)).\]
It follows $J^1_L(E)$ is a $(\B,\B)$-module and  a left $\B\otimes_{\A} \B$-module.
We get a canonical exact sequence of abelian groups
\begin{align}
&\label{Latiyah} 0\rightarrow E \rightarrow^i J^1_L(E) \rightarrow^j L\otimes_{\A} E\rightarrow 0
\end{align}
called the \emph{generalized Atiyah sequence associated to $L$}.
One checks that $i,j$ are maps of $(\B,\B)$-modules and hence maps of $\B\otimes_{\A} \B$-modules.

\begin{lemma} \label{LATsequence}
There is a one to one correspondence between the set of $(\B,\B)$-linear splittings of (\ref{Latiyah})
and the set of $L$-connections on $E$.
\end{lemma}
\begin{proof} Assume $s:L\otimes_{\A} E\rightarrow J^1_L(E)$ is a $(\B,\B)$-linear splitting of $j$.
It follows $s(x\otimes e)=(D(x\otimes e),x\otimes e)$. It follows
\[ D:L\otimes_{\A} E\rightarrow E\]
is a $(\B,\Der)$-linear map hence we get an $L$-connection $\nabla$. This gives a one to one correspondence
and the Lemma is proved.
\end{proof}

\begin{definition} The sequence (\ref{Latiyah}) gives rise to a cohomology class
\[a_L(E)\in \Ext^1_{\B\otimes_{\A} \B}(L\otimes_{\A} E,E)\]
called the \emph{generalized Atiyah class} of $E$.
\end{definition}

\begin{proposition} \label{aclass} The generalized Atiyah class $a_L(E)$ is zero if and only if $E$ has an $L$-connection
\[\nabla:L\rightarrow \End_{\A}(E).\]
\end{proposition}
\begin{proof}The proof follows from Lemma \ref{onetoone} and Lemma \ref{LATsequence}.
\end{proof}

Note: For any Lie-Rinehart algebra $L$ here is a canonical map of associative rings $p:\B\rightarrow \U(L)$.
A left $\B$-module $E$ has an $L$-connection 
\[ \nabla:L\rightarrow \End_{\A}(E) \]
if and only if $E$ is a left $\U(L)$-module. Hence for a left $\B$-module $E$ it follows 
the class $a_L(E)$ is zero if and only if the left $\B$-module structure on $E$
lifts to a left $\U(L)$-module structure on $E$. Hence the class $a_L(E)$ is the obstruction class for lifting
the left $\B$-module structure on $E$ to a left $\U(L)$-module structure on $E$.

\begin{example} The relative Kodaira-Spencer map and class.\end{example}

Let $E$ be an arbitrary $\B$-module and let $\mathbf{V}(E)\subseteq \Der_{\A}(\B)$ be the \emph{relative Kodaira-Spencer kernel}
of $E$ (see \cite{maa1}). By definition $\mathbf{V}(E)$ is the set of derivations $\delta \in \Der_{\A}(\B)$ with the property 
there is
a map $\nabla(x)\in \End_{\A}(E)$ such that 
\begin{align}
&\label{deriv}\nabla(x)(ae)=a\nabla(x)(e)+x(a)e
\end{align}

\begin{lemma} \label{lie} Let $x,y\in \Der_{\A}(\B)$ and let $\nabla(x), \nabla(y)\in \End_{\A}(E)$ 
be two maps satisfying equation
(\ref{deriv}). Let $\phi([x,y]):=[\nabla(x),\nabla(y)]$. It follows
\[ \phi([x,y])(ae)=a\phi([x,y])(e)+[x,y](a)e\]
for all $a\in \B$ and $e\in E$.
\end{lemma}
\begin{proof} The proof is an exercise.
\end{proof}

\begin{lemma} The subset $\mathbf{V}(E)\subseteq \Der_{\A}(\B)$ is a sub-Lie-Rinehart algebra of $\Der_{\A}(\B)$.
\end{lemma}
\begin{proof} One checks that $\mathbf{V}(E)$ is a sub-$\B$-module of $\Der_{\A}(\B)$. From Lemma \ref{lie} it follows
$\mathbf{V}(E)$ is a sub-$\A$-Lie algebra of $\Der_{\A}(\B)$. The Lemma follows.
\end{proof}
One may check that $\mathbf{V}(E)=ker(g)$ where
\begin{align}
&\label{KSmap} g:\Der_{\A}(\B)\rightarrow \Ext^1_{\B/\A}(E,E)
\end{align}
is defined as follows:
Let $x\in \Der_{\A}(\B)$. Let $E\oplus E$ have the following left $\B$-module structure:
\[ a(u,v):=(au+x(a)v,av).\]
Let $E(x)$ be the left $\B$-module $E\oplus E$ with $\B$-module structure defined by $x$.
It follows we get an exact sequence
\[ 0\rightarrow E \rightarrow E(x) \rightarrow E \rightarrow 0\]
of left $\B$-modules. Hence for any element $x\in \Der_{\A}(\B)$ we get an extension $g(x)\in \Ext^1_{\B/\A}(E,E)$.
It follows the map $g$ is a map of left $\B$-modules.
The map $g$ is the well known Kodaira-Spencer map from deformation theory. In \cite{maa1} the Kodaira-Spencer map
is defined using the Hochschild complex of $\Hom_k(E,E)$ where $k$ is a subfield of $\B$. The map from (\ref{KSmap})
is defined for an arbitrary map of commutative rings $\A\rightarrow \B$. It follows we may view the map
$g$ from (\ref{KSmap}) as the \emph{relative Kodaira-Spencer map}.

It follows we get a cohomology class
\begin{align}
&\label{ksclass}  a_{\mathbf{V}(E)}(E)\in \Ext^1_{\B\otimes_{\A} \B}(\mathbf{V}(E)\otimes_{\A} E, E).
\end{align}

\begin{corollary} Let $E$ be any $\B$-module and let $\mathbf{V}(E)$ be the relative Kodaira-Spencer kernel of $E$.
It follows $E$ has an $\mathbf{V}(E)$-connection if and only if $a_{\mathbf{V}(E)}(E)=0$.
\end{corollary}
\begin{proof} The proof follows from \ref{aclass}.
\end{proof}

Hence the class from (\ref{ksclass}) may be viewed as the \emph{relative Kodaira-Spencer class} of $E$. 
Note: In \cite{maa1} the Kodaira-Spencer class $a(E)$ is a class in the group 
\[ a(E)\in \Ext^1_B(\mathbf{V}(E), \End_{\B}(E)) ,\]
where $\B$ is a commutative ring over a field $k$. It is necessary for the proof of existence of the class in 
$\Ext^1_{\B}(\mathbf{V}(E), \End_{\B}(E))$ that $\B$ contains a field $k$. The class from (\ref{ksclass}) exists 
in greater generality. If we view the Kodaira-Spencer class as a class defined on $\B\otimes_{\A} \B$ we see it is a 
special case of the more general class $a_L(E)$ which is defined for an arbitrary Lie-Rinehart algebra $L$ and
an arbitrary $\B$-module $E$.

Dually there is for every morphism $\phi:\Omega^1_{\B/\A}\rightarrow \Omega$ of left $\B$-modules an Atiyah sequence
\begin{align}
&\label{atiyahclassical}0\rightarrow \Omega \otimes_{\B} E \rightarrow P^1_{\Omega}(E)\rightarrow E \rightarrow 0
\end{align}
defined using the first order module of principal parts $P^1_{\B/\A}(E)$ of $E$. Let
\[  D:\B\rightarrow \Omega \]
be defined by $D:=\phi \circ d$ where $d$ is the universal derivation. It follows $D$ is a derivation.
Define
\[ P^1_{\Omega}(E):=\Omega \otimes_{\B} E\oplus E \]
with the canonical left $\B$-module structure. Define the following right $\B$-module structure:
\[ (x\otimes e,f)b:=(x\otimes(eb)+D(b)\otimes f, fb).\]
It follows (\ref{atiyahclassical}) is right split by a $\B$-linear map $s:E\rightarrow P^1_{\Omega}(E)$ where
$s(e)=(\nabla(e),e)$ where $\nabla$ is an $\Omega$-connection.
An $\Omega$-connection is an $A$-linear map
\[ \nabla:E\rightarrow \Omega \otimes_{\B} E \]
where
\[ \nabla(ae)=a\nabla(e)+D(a)\otimes e.\]
The first order jet module $J^1_L(E)$ is in some cases the dual of the
first order module of principal parts $P^1_{\Omega}(E)$ of $\Omega$.

I belive the results in this example are well known but include them because of lack of a good reference.
The construction of sequence (\ref{atiyahclassical}) presented here follows the presentation in \cite{karoubi1}.

 Assume $D:L\otimes_A E\rightarrow E$ is an $(\B,\Der)$-linear map. We defined the \emph{curvature} 
of $D$ as follows:
\[ K_D:L\times L\times E\rightarrow E\]
\[K_D(x,y,e):=D([x,y]\otimes e)-D(x\otimes D(y\otimes e))+D(y\otimes D(x\otimes e)).\]
We say $D$ is \emph{flat} if $K_D=0$. It follows $D$ is flat if and only if the $L$-connection $\nabla$ associated
to $D$ is flat.
The curvature $K_D$ defines a  $B$-linear map
\[ K_D:L\wedge_{\B} L\otimes_{\B} E\rightarrow E.\]

\begin{example} \label{cohomology} \text{Cohomology of arbitrary $L$-connections.}

For any associative ring $R$ the category of left $R$-modules has enough injectives and one uses injective
resolutions of $R$-modules to construct functors $\Ext^i_R(M,N)$ for any pair of left $R$-modules $M,N$.

Let $\alpha:L\rightarrow \Der_{\A}(\B)$ be a Lie-Rinehart algebra and 
let $J\subseteq \U(L)$ be a 2-sided ideal. Let $\U_J(L)=\U(L)/J$.
It follows we may for any pair of left $\U_J(L)$-modules $V,W$ and any integer $i\geq 0$ define 
the $\Ext$-group $\Ext^i_{\U_J(L)}(V,W)$ as follows. Let 
\begin{align}
&\label{injective} 0\rightarrow W \rightarrow^{d_0} L^1 \rightarrow^{d_1} L^2 \rightarrow^{d_2} \cdots 
\end{align}
be an injective resolution of $W$ in the category of left $\U_J(L)$-modules. Apply the left exact functor $\Hom_{\U_J(L)}(V,-)$
to the sequence (\ref{injective}) to get the complex 
\begin{align}
&\label{cplx}0 \rightarrow \Hom_{\U_J(L)}(V,L^1)\rightarrow^{d_1^*} \Hom_{\U_J(L)}(V,L^2) \rightarrow^{d_2^*}  
\Hom_{\U_J(L)}(V,L^3) \rightarrow \cdots 
\end{align}
Note: The complex (\ref{cplx}) is by Lemma \ref{center} a complex of $\A$-modules.
\begin{definition} \label{extfunctor} Let $i\geq 0$ be an integer. Let
\[ \Ext^i_{\U_J(L)}(V,W)=ker(d_i^*)/im(d_{i-1}^*).\]
\end{definition}

Since $\conn(L)$ is a small abelian category, it follows from the 
Freyd-Mitchell full embedding theorem that there exists an associative  ring $R$ and an equivalence
between $\conn(L)$ and a sub-category of $\mod(R)$. 
The functor from the Freyd-Mitchell full embedding theorem does not preserve injective and projective objects hence 
we cannot in general use the ring $R$ from the theorem to define Ext and Tor functors for connections. 
Using the algebra $\U(L)$ we get an exact equivalence between
$\conn(L)$ and $\mod(\U(L))$. It follows projective and injective objects in $\conn(L)$ equals projective and injective
objects in $\mod(\U(L))$. Hence Definition (\ref{extfunctor}) is well defined.
It follows Definition (\ref{extfunctor}) gives the first definition of Ext-functors for 
arbitrary connections with no condition on the curvature.
\end{example}

\begin{example}\label{globalext}\text{ A globalization of Theorem \ref{equiv} and Definition \ref{extfunctor}.}

For any ringed topologiocal space $(X, \O_X)$ where $\O_X$ is an arbitrary sheaf of associative unital rings, 
it follows the category of $\O_X$-modules has enough injectives. Hence we may use the formalism of derived functors
to construct the derivative of any left exact functor $F:\mod(\O_X)\rightarrow \underline{C}$ where
$\mod(\O_X)$ is the category of left $\O_X$-modules and $\underline{C}$ is some abelian category.
There is work in progress constructing a globalization of Theorem \ref{equiv} to the case of an arbitrary
sheaf of Lie-Rinehart algebras $L$ on an arbitrary scheme $X$. One seek to define for any sheaf of Lie-Rinehart algebras $L$
on $X$ a sheaf of associative unital algebras $\mathcal{U}^{ua}_L$ on $X$ and an exact equivalence of categories
\[ F:\conn(L)\rightarrow \mod(\mathcal{U}_L^{ua}) \]
preserving injective objects. The space $(X, \mathcal{U}_L^{ua})$ will be a ringed space, and we may for any left exact
functor
\[ H: \mod(\mathcal{U}_L^{ua})\rightarrow \underline{C}, \]
where $\underline{C}$ is some abelian category, construct the derivatives of $H$. 
The aim is to define $\Ext$-functors for connections in the global situation. This is work in progress.
\end{example}

\begin{example}\label{motives}\text{ Some comments on Weil cohomology theories and motives.}

Let $\underline{Sm}_k$ be the category of smooth projective varieties over a field $k$ and let 
$\underline{Gr}(K)$ be the category of graded commutative $K$-algebras, where $K$ is a field.
A \emph{Weil cohomology theory} is a contravariant functor
\[ \H^*:\underline{Sm}_k\rightarrow \underline{Gr}(K) \]
satisfying a set of axioms including a K\"{u}nneth formula, Poincare duality and the Hard Lefschetz Theorem 
(see \cite{milne} for a precise description). Singular cohomology, DeRham cohomology, 
l-adic cohomology and cristalline cohomology are all examples of Weil cohomology theories.
One may ask if some these theories can be constructed using a global version of the $\Ext$-groups introduced in 
Defintion \ref{extfunctor}.
The theory of \emph{motives} is about constructing a \emph{universal Weil cohomology theory} ``parametrizing''
all Weil cohomology theories.
So far there is a construction of the derived category of motives. It is conjectured this category is equivalent
to the derived category of a Tannakian category of motives. 
\end{example}

By definition $\Ext^0_{\U_J(L)}(V,W):=\Hom_{\U_J(L)}(V,W)$.
It follows $\Ext^i_{\U_J(L)}(V,W)$ is an $\A$-module for every $i\geq 0$. 
The $\Ext$-groups are not left $\U_J(L)$-modules, hence there is in general no $L$-connection
\[ \nabla:L\rightarrow \End_{\A}(\Ext^i_{\U_J(L)}(V,W)).\]

By general results of \cite{cartan} we get for any short exact sequence of $L$-connections
\[ 0 \rightarrow (V', \nabla') \rightarrow (V,\nabla) \rightarrow (V'', \nabla'') \rightarrow 0\]
and any $L$-connection $(W, \rho)$, a long exact sequence of left $\A$-modules
\[ 0 \rightarrow \Hom_{\U(L)}(V'',W)\rightarrow \Hom_{\U(L)}(V,W)
\rightarrow \Hom_{\U(L)}(V',W)\rightarrow \]
\[ \Ext^1_{\U(L)}(V'', W)\rightarrow \Ext^1_{\U(L)}(V,W)
\rightarrow \Ext^1_{\U(L)}(V', W)\rightarrow \cdots .\]
Here we have chosen $J:=(0)$ to be the zero ideal.

In the relative setting there is a definition of $\Ext$-sheaves (see \cite{maa15}). 
The $\Ext$-sheaves defined in \cite{maa15} are equipped with a relative Gauss-Manin-connection.

Lie-Rinehart cohomology $\H^i(L,W)$ is in \cite{rinehart}  defined for every integer $i\geq 0$ as follows:
\[ \H^i(L,W):=\Ext^i_{U(\B,L)}(\B,W), \]
where $U(B,L)$ is Rineharts universal enveloping algebra of $L$ and $W$ is a flat $L$-connection.
When $L$ is a projective $\B$-module it follows $\H^i(L,W)$ may be calculated using the Lie-Rinehart complex.

\begin{lemma} Let $J=I_{L,0}$ and 
let $(W,\nabla)$ be a flat $L$-connection. It follows there is an isomorphism
\[ \Ext^i_{\U_J(L)}(\B,W)\cong \H^i(L,W)\]
for all integers $i\geq 0$.
\end{lemma}
\begin{proof} Since $\U_J(L)\cong U(\B,L)$ it follows
\[ \Ext^i_{\U_J(L)}(\B,W)\cong \Ext^i_{U(\B,L)}(\B,W)\cong \H^i(L,W).\]
The Lemma follows.
\end{proof}
Hence the group $\Ext^i_{\U_J(L)}(V,W)$ generalize Lie-Rinehart cohomology for all $i\geq 0$.
Using $\U_J(L)$ when we vary the ideal $J$, we get a cohomology theory $\Ext^i_{\U_J(L)}(V,W)$ defined for 
arbitrary $L$-connections $V,W$ with no conditions on the curvature. Definition (\ref{extfunctor}) 
gives a simultaneous construction of a large class of cohomology theories.

For a non-flat $L$-connection $\nabla:L\rightarrow \End_{\A}(W)$ we may use the Lie-Rinehart construction to get a sequence
of $\B$-modules
\begin{align}
&\label{lrcomplex} \cdots \rightarrow \Hom_{\B}(\wedge^p_{\B} L, W) \xrightarrow{d^p} \Hom_{\B}(\wedge^{p+1}_{\B}L, W)
\xrightarrow{d^{p+1}}\cdots
\end{align}
and maps of $\A$-modules. The sequence (\ref{lrcomplex}) is a complex if and only if $\nabla$ is a flat connection. Hence
for non-flat connections the sequence (\ref{lrcomplex}) does not give rise to well defined cohomology groups.

If the base ring $\A$ contrains a field $k$ we may consider the Hochschild cohomology groups of the 
left and right $\U_J(L)$-module $\Hom_k(V,W)$ where $V,W$ are left $\U_J(L)$-modules. We get an isomorphism
\begin{align}
\label{hochschild} \Ext^i_{\U_J(L)}(V,W)\cong \HH^i(\U_J(L), \Hom_k(V,W)) 
\end{align}
for all $i\geq 0$. It follows we may use the Hochschild complex to calculate the group $\Ext^i_{\U_J(L)}(V,W)$. 

\begin{example} \text{$\Ext$-groups, Hochschild cohomology and singular cohomology.}

Let $\B$ be a finitely generated regular commutative algebra over the complex numbers and let $X:=\Spec(\B)$.
Let $X_{\Cx}$ be the underlying complex algebraic manifold of $X$. Let $L:=\Der_{\Cx}(\B)$ be the module 
of derivations of $\B$. It follows there is for every $i\geq 0$ an isomorphism
\[ \H^i_{sing}(X_{\Cx}, \Cx)\cong \Ext^i_{U(\B,L)}(\B,\B).\]
The group $\Ext^i_{U(\B,L)}(\B,\B)$ may be described as the group of equivalence classes of exact sequences
of flat $L$-connections
\[ 0\rightarrow \B \rightarrow (V_1,\nabla_1)\rightarrow \cdots \rightarrow (V_i,\nabla_i)\rightarrow \B \rightarrow 0.\]
It may be this description of $\H^i_{sing}(X_{\Cx}, \Cx)$ will be helpful in the description 
of the image of the Chern character and cycle map in the smooth affine case. One seek to generalize Theorem \ref{main}
and give a descrpition of the set of cohomology classes $c\in \H^{2i}_{sing}(X_{\Cx}, \Cx)$ with the property 
that there is a finitely generated projective $\B$-module $E$ with $c_i(E)=c$.
This topic will be studied in a forthcoming paper on the subject.
\end{example}

\begin{example} \text{An extension of singular cohomology of an affine algebraic variety.}

Let $\alpha:L\rightarrow \Der_{R}(A)$ be a Lie-Rinehart algebra and let $f\in \Z^2(L,A)$ be a 2-cocycle.
We get for every integer $i\geq 0$ a non-trivial exact sequence of $\Ext$-groups
\[ 0\rightarrow \Ext_{U(A,L,f)}^i(A,A) \rightarrow \Ext_{\U(L)}^i(A,A)\rightarrow \H^i_f(L,A)\rightarrow 0.\]
Hence there is for every 2-cocycle $f$ and every integer $i\geq 0$ an inclusion of $R$-modules
\[ \Ext_{U(A,L,f)}^i(A,A) \subseteq \Ext_{\U(L)}^i(A,A).\]
Hence the cohomology group $\Ext_{\U(L)}^i(A,A)$ is \emph{universal} in the sense that it contains
the cohomology groups $\Ext_{U(A,L,f)}^i(A,A)$ for all $f\in \Z^2(L,A)$. There is in particular an inclusion
\[ \Ext^i_{U(A,L)}(A,A)\subseteq \Ext^i_{\U(L)}(A,A) \]
of $R$-modules, and by definition $\H^i(L,A):=\Ext^i_{U(A,L)}(A,A)$ is the Lie-Rinehart cohomology of $L$ with values
in $A$. This is because $U(A,L):=U(A,L,0)$.

\begin{definition} Let $\alpha:L\rightarrow \Der_R(A)$ be a Lie-Rinehart algebra and let $(W,\nabla)$ be an $L$-connections.
Let for any integer $i\geq 0$
\[ \H^i_{ua}(L,W):=\Ext_{\U(L)}(A,W) \]
be the \emph{universal cohomology} of the $L$-connection $(W,\nabla)$.
\end{definition}
From the dicsussion above it follows there is for every integer $i\geq 0$ an inclusion of $R$-modules
\[ \H^i(L,W)\subseteq \H^i_{ua}(L,W), \]
hence the universal cohomology of the $L$-connection $(W,\nabla)$ contains the Lie-Rinehart cohomology group
$\H^i(L,W)$ in the case when $\nabla$ is flat.

Let $X:=\Spec(A)$ where $A$ is a finitely generated regular algebra over the complex numbers and let
$X_{\Cx}$ be the underlying complex algebraic manifold of $X$ in the strong topology. It follows we get
for all integers $i\geq 0$ an exact sequence of $\Cx$-vector spaces
\begin{align}
&\label{ext} 0 \rightarrow \H^i_{sing}(X_{\Cx}, \Cx) \rightarrow \H^i_{ua}(\Der_{\Cx}(A),A) 
\rightarrow \H^i_0(\Der_{\Cx}(A),A)\rightarrow 0, 
\end{align}

where $\U(\Der_{\Cx}(A))$ is the universal algebra of $\Der_{\Cx}(A)$. 
One may ask for a topological interpretation of the cohomology groups 
\[ \H^i_{ua}(\Der_{\Cx}(A), A)\text{ and }\H^i_0(\Der_{\Cx}(A),A) \]
for $i\geq 0$. 
\end{example}

We may similarly define the $\Tor$-functors $\Tor^i_{\U_J(L)}(V,W)$ for arbitrary $\U_J(L)$-modules $V,W$. For 
a treatment of general properties of $\Ext$ and $\Tor$-functors for modules on arbitrary associative rings 
see \cite{cartan}. If we let for any flat $L$-connection $(W,\nabla)$
\[ \H_i(L,W):=\Tor^i_{U(\B,L)}(\B,W)\]
we get Rineharts \emph{i'th homology group} of the $L$-connection $W$.
We get a definition of cohomology and homology for connections in complete generality using
the algebra $\U_J(L)$. 

In the case when the Freyd-Mitchell 
full embedding theorem gives an equivalenced of categories $\mod(R)\cong \conn(L)$
it follows we may use this equivalence to prove existence of cohomology and homology of arbitrary connections.
The proof of the theorem does not give a practical method to calculate a ring $R$ with $\mod(R)\cong \conn(L)$. 
Given an associative ring $R$ and an equivalence of categories
\[ \mod(R)\cong \conn(L)\]
it follows the ring $R$ is morita equivalent to the following ring:
\[ \Mat_{n_1}(\cdots \Mat_{n_k}(R)\cdots ) \]
where $n_i\geq 1$ are integers for $i=1,..,k$. Here $\Mat_n(R)$ is the ring of $n\times n$-matrices with coefficients in $R$. 
We get an equivalence of categories
\[ \mod(\Mat_{n_1}(\cdots \Mat_{n_k}(R)\cdots )) \cong \mod(L).\]
Hence the ring appearing in the proof of the Freyd-Mitchell full embedding theorem may be large.

Using the universal algebra $\U(L)$ and the quotient $\U_J(L)$ we get an elementary and explicit construction of such a ring. 
Because of Formula (\ref{hochschild}) it may be possible to explicitly calculate $\Ext^i_{\U_J(L)}(V,W)$ 
using Hochschild cohomology.

\begin{example} \text{The Cartan-Eilenberg complex of a connection.}

Let $\nabla:L\rightarrow \End_{\A}(W)$ be an $L$-connection and let
\[ D:L\otimes_{\A} E \rightarrow E \]
be the corresponding $(\B,\Der)$-linear map.
Recall that
\[ K_D(x\wedge y\otimes e):=D(x\otimes D(y\otimes e))-D(y\otimes D(x\otimes e))-D([x,y]\otimes e)\]
is the curvature of $D$.
Define the following map
\begin{align}
&\label{koszul} d_p:\wedge_{\A}^{p} L\otimes_{\A} E\rightarrow \wedge^{p-1}_{\A}L\otimes_{\A} E 
\end{align}
for $p\geq 1$ an integer, by
\[ d_p(x_1\wedge \cdots \wedge x_p\otimes e):=\]
\[\sum_{i=1}^p(-1)^{i+1}x_1\wedge \cdots \wedge\overline{x_i}\wedge \cdots \wedge x_p\otimes D(x_i\otimes e)+\]
\[\sum_{i<j}(-1)^{i+j}[x_i,x_j]\wedge x_1\wedge \cdots \wedge \overline{x_i}\wedge \cdots \wedge
\overline{x_j}\wedge \cdots \wedge x_p\otimes e.\]
Let $d_0=D$ where $D:L\otimes_{\A} E\rightarrow E$.
The following formula is classical:
\[ d_{p-1}\circ d_p(x_1\wedge \cdots \wedge x_p\otimes e)=\]
\[ \sum_{i<j}(-1)^{i+j-1}\wedge x_1\wedge \cdots \wedge \overline{x_i}\wedge \cdots \wedge
\overline{x_j}\wedge \cdots \wedge x_p\otimes K_D(x_i\wedge x_j\otimes e).\]
Hence the maps from (\ref{koszul}) defines a complex if and only if $K_D=0$. If $L$ is a projective $\B$-module it follows
from \cite{rinehart} the homology of (\ref{koszul}) calculates homology $\H_i(L,W)=\Tor_{U(\B,L)}^i(\B,W)$ 
of the flat connection $W$. If $W$ is not flat the complex (\ref{koszul}) does not give rise to well defined homology groups. 
If $\A$ contains a field $k$, it follows Hochschild homology of the left and right $\U_J(L)$-module $\Hom_k(V,W)$ calculates
the group $\Tor^i_{\U_J(L)}(V,W)$ for any $\U_J(L)$-modules $V,W$.
If $\B=\A=k$ is a field, $L$ a $k$-Lie algebra and $E$ a left $L$-module, it follows 
(\ref{koszul}) is the classical Cartan-Eilenberg complex of $E$ computing the homology $\H_i(L,E)$.
\end{example}

\begin{corollary} Let $J=I_{L,0}$ and let $W$ be a flat $L$-connection. It follows
\[ \Tor^i_{\U_J(L)}(\B,W)\cong \H_i(L,W) \]
for all integers $i\geq 0$.
\end{corollary}
\begin{proof} Since $\U(L)/J\cong U(\B,L)$ the result follows by definition of $\H_i(L,W)$.
\end{proof}

Hence the $\A$-modules $\Ext^i_{\U_J(L)}(V,W)$ and $\Tor^i_{\U_J(L)}(V,W)$ generalize Lie-Rinehart cohomology and homology
as defined in \cite{rinehart}.

It might be Corollary \ref{equiv}
will be helpful in giving explicit calculations of the Grothendieck group $\K_0(\conn^{fp}(L))$. The group 
$\K_0(\conn^{fp}(L))$ is an important invariant for the Lie-Rinehart algebra $L$ and it is an open problem to 
calculate it explicitly.

\begin{example} \text{Moduli spaces of connections.}

In \cite{simpson} Simpson constructs moduli spaces of sheaves of modules on sheaves of rings of differential operators
on smooth projective varieties over the complex numbers. In the affine situation one may try to construct
moduli spaces of left $U(\B,L,f)$-modules and moduli spaces of left $\U(L)$-modules. Since a left $U(\B,L,f)$-module corresponds
to a left $\U(L)$-module it follows the moduli space of left $\U(L)$-modules ``contains'' the moduli space of left
$U(\B,L,f)$-modules. If one globalize this situation and constructs for any finite rank locally free sheaf of 
Lie-Rinehart algebras $L$ on a projective scheme $X$ over a field, the moduli space of coherent sheaves of left $\U(L)$-modules, one gets a generalization of Simpson's moduli space. Many problems in the theory of connections reduce to proving existence 
of connections where the curvature of the connection satisfy various properties. Hence such a construction 
might be useful in applications of the theory.
\end{example}

\begin{example} \text{Grothendieck groups of categories of $L$-connections.}

Recall the following construction from the previous chapter:
\[ \overline{\K}(L)_{\mathbf{Q}}=\{ \sum_i r_i[E_i,\nabla_i]: r_i\in \mathbf{Q}, R_{\nabla_i}=f_i, f_i\in \Z^2(L,\B)\}.\]
By Corollary \ref{equiv} there is an isomorphism of rings
\[ \K(L)_{\mathbf{Q}}\cong \K(\mod^{fp}(\U(L)))_{\mathbf{Q}} \]
and an inclusion of rings
\begin{align}
&\label{cherncharacter} \overline{\K}(L)_{\mathbf{Q}} \subseteq \K(\mod^{fp}(\U(L)))_{\mathbf{Q}}.
\end{align}
If we can give generators for $\K(\mod^{fp}(\U(L)))_{\mathbf{Q}}$ as left $\overline{\K}(L)_{\mathbf{Q}}$-module
me may use this description to study the image of the Chern character $Ch$ in $\H^*(L,\B)$ as discussed
in the previous section of this paper.
\end{example}

\end{document}